\documentclass[12pt]{article}

\usepackage[T1]{fontenc}

\usepackage{comment}

\usepackage{graphicx} 
\usepackage{amsthm,amsmath,amsfonts,amssymb,mathrsfs}
\usepackage{hyperref}
\usepackage{bbm}
\usepackage[letterpaper,top=2cm,bottom=2cm,left=3cm,right=3cm,marginparwidth=1.75cm]{geometry}
\usepackage{tabularx}

\usepackage{caption}
\usepackage{subcaption}

\usepackage{bbold} 
\usepackage[normalem]{ulem} 
\usepackage{placeins} 

\newtheorem{theorem}{Theorem}
\newtheorem{lmm}[theorem]{Lemma}
\newtheorem{prpstn}{Proposition}

\newtheorem{dfntn}[theorem]{Definition}

\newtheorem{rmrk}{Remark}

\newcommand{\xLone}{L^1}
\newcommand{\xLtwo}{L^2}
\newcommand{\xHone}{H^1}
\newcommand{\xHn}[1]{H^{#1}}
\newcommand{\xCone}{C^1}
\newcommand{\xCn}[1]{C^{#1}}

\newcommand{\Esp}{\mathbb{E}}
\newcommand{\Var}{\mathbb{V}\mathrm{ar}}
\newcommand{\tot}{\textrm{tot}}
\newcommand{\R}{\mathbb{R}}
\newcommand{\N}{\mathbb{N}}
\newcommand{\domain}{\mathbb{X}}
\newcommand{\probClass}{\mathrm{P}(a,b)}
\newcommand{\weight}{w}
\newcommand{\wlin}{\weight_{\mathrm{lin}}}
\newcommand{\LoneLeb}{\xLone(a,b)}
\newcommand{\Ltwo}{\xLtwo(\mu)}
\newcommand{\LtwoW}{\xLtwo(\mu, \weight)}
\newcommand{\HoneW}{\xHone(\mu, \weight)}
\newcommand{\Hm}[1]{\xHn{{#1}}(\mu)}
\newcommand{\Cm}[1]{\xCn{{#1}}}
\newcommand{\LtwoDot}[2]{\langle #1, #2 \rangle }
\newcommand{\LtwoDotW}[2]{\langle #1, #2 \rangle_{\weight} }
\newcommand{\LtwoDotH}[2]{\langle #1, #2 \rangle_{\HoneW} }
\newcommand{\LtwoNorm}[1]{\Vert #1 \Vert}
\newcommand{\LtwoNormW}[1]{\Vert #1 \Vert_{\weight} }
\newcommand{\HoneNorm}[1]{\Vert #1 \Vert_{\HoneW}}
\newcommand{\funClassWeight}{\mathcal{\mathrm{W}}(a,b)}
\newcommand{\eigfun}{\varphi}
\newcommand{\eigval}{\lambda}
\newcommand{\indBasis}{j}  
\newcommand{\indVar}{k}  
\newcommand{\indDoE}{i}  
\newcommand{\p}{p} 
\newcommand{\enum}{, \ldots ,}
\newcommand{\ve}[1]{\boldsymbol{#1}}
\newcommand{\ca}{\mathcal{A}}
\newcommand{\cm}{\mathcal{M}}
\newcommand{\cu}{\mathcal{U}}
\newcommand{\cx}{\mathcal{X}}
\newcommand{\innprod}[3][]{\left\langle #2, #3 \right\rangle_{#1}}
\newcommand{\partderiv}[2][]{\frac{\partial #1}{\partial x_{#2}}}
\newcommand{\norme}[2]{\left\| #1 \right\|_{#2}}

\renewcommand{\rho}{r}

\usepackage{xcolor}

\newcommand{\alp}{{\ve \alpha}}

\newcommand{\rev}[1]{{\color{black}{#1}}}

\usepackage{authblk} 

\title{Gradient-enhanced global sensitivity analysis\\ with Poincar\'e chaos expansions}

\author[1]{O. Roustant}
\author[2]{N. L\"uthen}
\author[1]{D. Heredia}
\author[2]{B. Sudret}

\affil[1]{{\small UMR CNRS 5219, Institut de Math\'ematiques de Toulouse, INSA, Universit\'e de Toulouse, France}}
\affil[2]{\small Chair of Risk, Safety and Uncertainty Quantification, ETH Z\"urich, 8093 Z\"urich, Switzerland}

\date{\today}

\begin{document}

\maketitle

\begin{abstract} Spectral methods, also known as chaos expansions, are widely used in global sensitivity analysis (GSA), as they leverage orthogonal bases of $\xLtwo$ spaces to efficiently compute Sobol' indices, particularly in data-scarce settings. When derivatives of the model are available, a desirable property,
\rev{both for modeling and GSA purposes,} is for the derivatives of the basis functions to also form an orthogonal basis. 
We demonstrate that the only basis satisfying this property is the one associated with weighted Poincaré inequalities and Sturm–Liouville eigenvalue problems, which we call 
Poincaré basis. \rev{We also show that under certain conditions the Poincaré basis achieves the same convergence rate as the best polynomial approximation for classes of smooth functions.}

We then introduce a comprehensive framework for gradient-enhanced GSA 
that integrates recent advances 
both in the construction of the expansion - with gradient-enhanced regression - and in the construction of weights for derivative-based sensitivity analysis. Furthermore, the proposed methodology is applicable to a broad class of probability measures and various choices of weights. 
We illustrate its efficiency on a challenging flood modeling case study, where Sobol' indices are accurately estimated using limited data. 
\end{abstract}

\section{Introduction}

The analysis of complex input/output systems has received growing attention in the last decades. 
Here we will consider a system that can be described as a multivariate real-valued function $\cm: \ve{x} \in \domain \subset \R^d \to \R$. Global sensitivity analysis (GSA) aims at quantifying the influence of some input variable $X_i$, viewed as a random variable, on the variability of the output $\cm(\ve{X})$ (with $\ve{X} = (X_1, \dots, X_d)$). 
\rev{We further assume that $X_1, \dots, X_d$ are independent and we denote by $\mu$ the product probability measure of $\ve{X}$}.
Famous indicators are variance-based sensitivity indices, also called Sobol' indices. Beyond their simplicity, one reason for their success is the existence of the Sobol-Hoeffding decomposition of $\cm(\ve{X})$ as a sum of \emph{orthogonal} terms corresponding to main effects and interactions. Chaos expansion methods, which rely on multivariate orthonormal bases, are particularly suitable to compute Sobol' indices, as they can intrinsically leverage this orthogonality property, as originally shown in \cite{SudretRESS2008}.
\rev{We note in passing that it is not recommended to use the quantile transformation to get uniform random variables $U_i = F_i^{-1}(X_i)$ (where $F_i$ is the cumulative density function of $X_i$) and to consider the new model $\tilde{\mathcal{M}}(\mathbf{U}) = \mathcal{M}(F_1(U_1), \dots, F_d(U_d))$. Such a transformation would induce strong non-linearities, inherited from the $F_i$'s even if the original model $\mathcal{M}$ is linear, \rev{and may make the model more difficult to approximate (see, e.g., \cite{torre2019data}, \S 2.3.2.)}. Thus, it is \rev{preferable} to construct orthonormal bases in $\Ltwo$ when $\mu$ is not uniform, and go beyond some common choices such as Fourier or wavelet bases.}

From now on, we assume that the gradient of $\cm$ is available everywhere on $\domain$, and we aim at using this gradient to improve the computation of Sobol' indices with chaos expansions. 
In general, chaos expansion methods rely on an orthonormal basis $(\psi_\indBasis)_{\indBasis \in \N}$ of $\Ltwo$, built by tensorization of univariate orthonormal bases associated to each input variable.
Then any $\cm$ in $\Ltwo$ can be expanded as
\begin{equation} \label{eq:chaosExpansion}
   \cm(\ve{x}) = \sum_{\indBasis \in \N} c_\indBasis \psi_\indBasis(\ve{x}).
\end{equation}
When the gradient of $\cm$ is available, a desirable property is that for all $\indVar=1, \dots, d$, the partial derivative\rev{s $\left(
\frac{\partial \psi_j}{\partial x_\indVar}\right)_{j \in \mathbb{N}}
$ also form} an orthogonal basis. 
This is due to two main reasons. 
First, this enhances the estimation of the expansion coefficients by sparse regression methods using both function and derivative values. In this setting \cite{Adcock_2019} obtains theoretical recovery guarantees, by using a suitable notion of coherence.
Second, this is also a favorable situation for GSA, where Sobol' indices are computed with Parseval's formula from the basis expansion \eqref{eq:chaosExpansion}. Indeed, in that case a derivative of the expansion \eqref{eq:chaosExpansion} also appears as a basis expansion
\begin{equation} \label{eq:chaosExpansionDer}
    \frac{\partial \cm}{\partial x_i}(\rev{\ve x}) = \sum_{\indBasis \in \N} c_\indBasis \frac{\partial \psi_j}{\partial x_i}(\rev{\ve x}), 
\end{equation}
from which Sobol' indices can be derived \cite{PoincareChaos}. When the derivative varies less than the function, this leads to more accurate estimations \cite{PoincareChaos, Poincare_Chaos_Sparse_2023}.

As argued in the last paragraph, we are looking for orthogonal bases of $\Ltwo$ that are stable by derivation, in the sense that their partial derivatives also form an orthogonal basis.  
By construction of the chaos expansion, it is sufficient to restrict this problem to univariate bases. It is also convenient to add a degree of freedom by considering a different Hilbert space for the derivatives, such as the weighted $\xLtwo$ Hilbert space
$ \LtwoW = \{f: \domain \to \R, \text{measurable, s.t. } \int_{\domain} f^2(\ve{x}) w(\ve{x}) \, d\mu(\ve{x}) < +\infty\} $, where $w$ is some positive function. \rev{Note that, apart from specific cases (recalled below) the bases formed by orthogonal polynomials do not satisfy this property.}
In this paper, we prove that the only one-dimensional orthonormal basis of $\Ltwo$ such that its derivative is an orthogonal basis of $\LtwoW$ is formed by eigenfunctions of the spectral problem associated to weighted Poincaré inequalities. This generalizes a previous result, in the case of unweighted Poincaré inequalities \cite{Poincare_Chaos_Sparse_2023}. Such a basis, called Poincaré basis, also corresponds to the eigenfunctions of a Sturm-Liouville eigenvalue problem, and was studied in \cite{Adcock_2019}. We call Poincaré chaos expansion the chaos expansion method associated to the Poincaré basis.

\rev{While our construction of the Poincaré basis is guided by gradient-based sensitivity analysis, it is equally important to assess its approximation properties, and investigate the convergence rate of the basis expansion for a given class of functions. These properties are derived from Sturm-Liouville theory. When the probability distribution $\mu$ has finite support with probability density function $r$, we show that
whenever the product $wr$ vanishes at the boundaries, 
the Poincaré basis 
achieves the same convergence rate as the best polynomial approximation, for classes of functions with a given level of smoothness.
This insight provides practical guidance for selecting $w$: choosing a weight function that vanishes at the boundary ensures strong approximation properties for the corresponding Poincaré basis, regardless of $\mu$.}

In general, the Poincaré basis does not coincide with orthogonal polynomials, except for three particular cases corresponding to Hermite, Laguerre and Jacobi polynomials (the latter includes Legendre and Chebyschev as special cases). Interestingly, these cases are the ones traditionally considered in gradient-enhanced polynomial chaos expansions \cite{Jakeman2015, Hampton2015, Peng2016, Guo2018}. Thus, a common point between them is the stability property of the orthogonal basis with respect to derivation, which enhances both the basis expansion and its usage in GSA, as explained above. Hence, considering Poincaré chaos expansions gives a wider natural framework for leveraging gradient information.

In the second part of the paper, we develop a comprehensive framework that integrates recent work on the computation of Poincaré chaos expansions by sparse methods with the construction of weighting schemes for derivative-based sensitivity analysis. As observed in \cite{HerediaWeightPoincare}, choosing the weight $w$ of $\LtwoW$ can be beneficial from a GSA perspective. 
Concerning the computation by sparse methods, we consider both the multi-output regression of \cite{Adcock_2019}, involving a large regression matrix, and the aggregation of multiple single-output regressions, with smaller regression matrices, as proposed in \cite{Poincare_Chaos_Sparse_2023}. 
To promote its usage among researchers and practitioners, 
the whole methodology is implemented in an open-source software. The numerical methods used allow us to go beyond the known analytical cases. Thus, we can deal with a broad class of probability measures and various weight choices. 
We put in action the whole methodology in a challenging case study, where some input variables follow unusual truncated probability distributions.

The paper is organized as follows. Section~\ref{sec:background} reviews general concepts related to chaos expansions and global sensitivity analysis. Section~\ref{sec:poincare_basis} introduces the univariate Poincaré basis, which is associated with weighted Poincaré inequalities. It provides sufficient conditions for its existence and presents a key characterization of this basis among orthogonal bases: its stability under differentiation. 
\rev{Under additional conditions, we also
address the approximation properties of this basis.}
Section~\ref{sec:PoinCE} focuses on chaos expansions constructed from Poincaré bases. It presents a closed-form expression for derivative-based sensitivity measures, along with two methods for estimating the chaos coefficients. Section~\ref{sec:numExp} illustrates the performance of the Poincaré chaos expansion in a gradient-enhanced setting through numerical experiments, with particular attention to a flood risk case study involving various non-standard probability distributions.
For readability, \rev{technical details and} all proofs are postponed to the appendix.

\section{Background}
\label{sec:background}

\rev{
\subsection{Generalized chaos expansions}
\label{sec:CEgeneral}
}
We consider a model $\cm : \R^d \to \R$ which has finite variance under the joint probability \rev{measure} $\mu$ of the input random variables $\ve X$ --- in other words, $\cm \in \Ltwo$. We assume the input random variables to be independent, therefore $\mu$ \rev{is a product measure:  $\mu = \otimes_{k = 1}^d \mu_\indVar$}.

Let $\{ \psi_\indBasis: \indBasis \in \N \}$ be an orthonormal basis of $\Ltwo$, i.e., a complete orthonormal system satisfying
\[
    \innprod{\psi_i}{\psi_\indBasis} = \int \psi_i(\ve x) \psi_\indBasis (\ve x) \, d\mu (\ve x)= \begin{cases}
        0 \text{ if } i \neq j, \\
        1 \text{ else.}
    \end{cases}
\]
Then any $\cm \in \Ltwo$ can be expanded as
\[
    \cm(\ve x) = \sum_{\indBasis \in \N} \innprod{\cm}{\psi_\indBasis} \psi_\indBasis(\ve x).
\]
In the context of uncertainty quantification, this is usually called a \rev{(generalized)} \emph{chaos expansion} \rev{\cite{Wiener1938,Ghanembook1991,Ernst2012,PoincareChaos}}.

A standard choice of chaos expansion is the \emph{polynomial chaos expansion} (PCE) which makes use of orthonormal polynomials \cite{XiuKarniadakis2002,Soize2004}. Let $\{\psi_{\indVar,\indBasis}, \indBasis \in \N\}$ denote univariate polynomial basis which are orthonormal with respect to $\mu_\indVar$, $\indVar = 1 \enum d$, and $\indBasis$ denotes the degree. For all univariate polynomial bases, we set $\psi_{\indVar,0} = 1$. The multivariate orthonormal basis is constructed from the univariate bases as the following tensor product:
\begin{equation} \label{eq:CEdef}
    \psi_\alp(\ve x) = \prod_{\indVar=1}^d \psi_{\indVar, \alpha_\indVar}(x_\indVar),
\end{equation}
where $\alp \in \N^d$ is called a \emph{multi-index} and characterizes the degree of the basis polynomial in each of the input variables. The \emph{total degree} of a basis polynomial is defined by $\sum_{\indVar=1}^d \alpha_\indVar$. 

A second choice of chaos expansions, which has been proposed in \cite{PoincareChaos} and explored further in \cite{Poincare_Chaos_Sparse_2023} and \cite{HerediaWeightPoincare}, are the \emph{Poincar\'e chaos expansions} (PoinCE).
Similarly to PCE, they are constructed by tensorization of univariate orthonormal bases with Eq.\eqref{eq:CEdef}. The computation of the univariate bases and the resulting properties are explained in Section~\ref{sec:poincare_basis}.
Poincar\'e chaos expansions were studied independently by \cite{Adcock_2019}, who consider mainly the special case when the basis is polynomial.

\subsection{Variance-based and derivative-based SA}
\label{sec:GSAreminder}
In this section, we recall general concepts from GSA, referring to \cite{GSAbook} for more details.

\emph{Variance-based sensitivity indices.}
Variance is one of the simplest indicators of variability, and variance-based sensitivity indices, known as Sobol' indices, are logically the first quantities of interest in GSA. They rely on the Sobol'-Hoeffding decomposition of $\cm \in \Ltwo$ under the assumption of independent input variables, written as
$$ \cm(\ve{X}) = \sum_{I \subseteq \{1, \dots, d \}} \cm_I(\ve{X}_I).
$$
Here $X_I$ denotes the sub-vector of $X$ obtained by selecting the coordinates that belong to $I =\{i_1, \, i_2, \cdots,i_s\} \subseteq \{1, \dots, d \}$. The decomposition is unique under the non-overlapping condition $\Esp(\cm_I(\ve{X}_I) \vert X_J)=0$ for all strict subsets $J \subset I$, with the convention $\Esp(. \vert X_\emptyset) = E(.)$. In that case, the terms are orthogonal, which allows to decompose the variance of $\cm(\ve{X})$ as a sum of components associated to sets of variables. Variance-based sensitivity indices are then defined as ratios of variance. For a single variable $X_\indVar$, the Sobol' index $S_\indVar$ and the total Sobol' index $S_\indVar^\tot$ are defined by
$$  S_\indVar = \frac{\Var \left[ \cm_\indVar(\ve{X}) \right]}{\Var \left[ \cm(\ve{X}) \right]}, \qquad 
S_\indVar^\tot = \frac{ \sum_{I \supseteq \{ \indVar \} }\Var \left[ \cm_I(\ve{X}) \right]}{\Var \left[ \cm(\ve{X}) \right]}.$$
Here we will focus on the total Sobol' index, which can be used to detect inactive (also called unimportant) variables. Indeed, under mild conditions on $\cm$ and $\mu$, if $S_\indVar^\tot = 0$ then $\cm$ does not depend on $x_\indVar$.

\emph{Derivative-based sensitivity measures.}
When the gradient of $\cm$ is provided, global sensitivity indices can be obtained by integration of local ones based on partial derivatives. We will consider here the (weighted) derivative-based sensitivity measure (DGSM, \rev{\cite{SobolKucherenko2009}}) associated to a single variable $X_\indVar$, of the form
$$ \nu_\indVar = \Esp \left[ w_\indVar(X_\indVar) \left(  \partderiv[\cm]{\indVar}(\ve{X}) \right)^2 \right],$$
where $w_\indVar$ is a non-negative function. Similarly to total Sobol' indices, DGSM can be used to detect inactive variables: under mild conditions on $\cm$ and $\mu$, if $\nu_\indVar = 0$, then $\cm$ does not depend on $x_\indVar$.

\emph{Link with chaos expansions.} Chaos expansions are particularly suitable to compute Sobol' indices, as they can leverage the orthogonality of the Sobol' decomposition. Thus, once the chaos expansion $\cm(\ve x) = \sum_{\alp \in \N^d} c_\alp \psi_\alp(\ve x) $ has been computed -- typically by using a sparse regression technique (see Section~\ref{sec:sparse_regression}), all Sobol' indices are computed from mere sums of squared coefficients \cite{SudretRESS2008}. For instance,
\begin{equation} \label{eq:totalSobolFromCE}
S_\indVar^\tot = \frac{\sum_{\alp \in \N^d, \alp_\indVar \geq 1} c_\alp^2}{\sum_{\alp \neq \ve 0} c_\alp^2}.
\end{equation}
In practice, estimates of $S_\indVar^\tot$ are obtained by restricting the summations to terms of the truncated expansion defined by a \rev{finite subset of multi-indices} $\mathcal{A}$. 
On the other hand, since the partial derivatives $\partderiv[\psi_\alp]{\indVar}$ do not form an orthogonal basis in general, DGSM may not be simply derived from chaos expansions. \rev{Note that this can be achieved by projecting the derivative basis functions onto the original basis and recombining the PCE coefficients, as shown in \cite{SudretMaiRESS2015}.
However, this approach involves rather cumbersome bookkeeping, and closed-form formulas are available only in a few cases.}

\section{Poincar\'e basis on the real line}
\label{sec:poincare_basis}

\subsection{Setting and notations}
Let $(a,b)$ be an open interval of the real line, with $-\infty \leq a < b \leq +\infty$. When $a$ and/or $b$ are infinite, we adopt the convention that $[-\infty,b]=(-\infty,b]$ and/or $[a,\infty]=[a,\infty)$.\\
We will denote by $\LoneLeb$ the space of measurable functions $f:(a,b) \to \R$ which are integrable with respect to the Lebesgue measure: $\int_a^b \vert f(t) \vert dt < +\infty$.\\
For simplicity, in the sequel we will remove the integration variable in the integrals, denoting $\int_a^b f$ instead of $\int_a^b f(t) dt$ and $\int f d\mu$ instead of $\int f(t) d\mu(t)$.
\\
Following \cite{HerediaWeightPoincare}, we define:

\begin{itemize}
    \item $\probClass$: the set of probability measures $\mu$ on $(a,b)$ whose \rev{probability density function (pdf) with respect to the Lebesgue measure, denoted by} $\rho$, is continuous and piecewise $\xCone$ on $[a,b]$, positive on $(a,b)$. 
    \item $\funClassWeight$: the set of continuous functions on $[a,b]$, that are piecewise $\xCone$ and positive on $(a,b)$.
\end{itemize}

Notice that, compared to \cite{HerediaWeightPoincare}, we have slightly relaxed the assumptions on $\rho$ by allowing it to vanish at $a, b$.\\
For $\mu \in \probClass$ and $w \in \funClassWeight$, we consider the weighted $\xLtwo$ space
$$ \LtwoW = \{f: (a,b) \to \R, \text{measurable, s.t. } \int_a^b f^2 w \, d\mu < +\infty\} $$
with inner product $\LtwoDotW{f}{g} = \int f g \, \weight \, d\mu$. 
The usual (unweighted) $\xLtwo$ space is denoted by $\Ltwo \equiv \xLtwo(\mu, 1)$, with inner product $\langle., .\rangle$. 
Finally, we consider the weighted Sobolev space
$$
\HoneW=\left\{f\in \Ltwo, \textrm{ s.t. } f'\in \LtwoW\right\} $$
with inner product $\LtwoDotH{f}{g} = \LtwoDot{f}{g} + \LtwoDotW{f'}{g'}$. Here $f'$ stands for the weak derivative of $f$. As $\weight \, \rho$ is positive almost everywhere (with respect to the Lebesgue measure), $\LtwoW$ and $\HoneW$ are Hilbert spaces (see e.g. \cite{kufner_weighted_sobolev_spaces}).\\
The norms of $\Ltwo, \LtwoW, \HoneW$ are denoted respectively $\LtwoNorm{.}, \LtwoNormW{.}, \HoneNorm{.}$.

\subsection{Definition and characterization}
Let $\mu \in \probClass$ and $\weight \in \funClassWeight$. We say that $\mu$ satisfies a Poincar\'e inequality with weight $\weight$ if there exists a constant $C$ such that for every function $f\in \HoneW$ verifying $\int_a^b f \, d\mu =0$, we have
\begin{equation}
	\label{eq:PoincareInequality}
		\int_{a}^b f^2 \, d\mu 
		\leq C \int_a^b \weight \, (f')^2 \, d\mu.
\end{equation}

Poincaré inequalities are closely linked to the spectral problem of finding $\eigval \in \R$ and $f \in \HoneW$ such that
\begin{equation} \label{eq:PoincareWeakProblem}
	\LtwoDotW{f'}{g'} = \eigval \LtwoDot{f}{g}, \qquad \forall g \in \HoneW,
\end{equation} 
where the two $\xLtwo$ norms involved in Eq.\eqref{eq:PoincareInequality} have been replaced by their associated bilinear forms. Indeed, the smallest constant $C$ in Eq.\eqref{eq:PoincareInequality} (called Poincaré constant) corresponds to the inverse of the first non-zero eigenvalue in Eq.\eqref{eq:PoincareWeakProblem}, when these quantities exist
(see \textit{e.g.} \cite{BGL}).
In our context, we are interested in the eigenfunctions of Eq.\eqref{eq:PoincareWeakProblem}, that we call Poincaré basis, extending the definition used in \cite{Poincare_Chaos_Sparse_2023} in the case $\weight \equiv 1$. 

\begin{dfntn}[Poincaré basis] \label{def:PoincareBasis}
Consider the spectral problem \eqref{eq:PoincareWeakProblem}.
Assume that there exists a countable set of eigenvalues
$(\eigval_\indBasis)_{\indBasis \in \N}$, with ${0 = \eigval_0 < \eigval_1 < \cdots}$ and a countable set of eigenfunctions $(\eigfun_\indBasis)_{\indBasis \in \N}$ constituting an orthonormal basis of $\Ltwo$.
Then we call that basis of eigenfunctions \emph{Poincaré basis}.
\end{dfntn}

Notice that as the eigenvalues are all simple, the Poincaré basis is uniquely defined, up to a change \rev{of} sign of its elements.\\
Poincaré inequalities are closely related to the diffusion operator
\rev{$L_w\colon D(L_w) \subseteq \Ltwo \to \Ltwo$ defined by
\begin{equation} \label{eq:DiffusionOperator}
	L_w(f) = \frac{1}{\rho} (\weight \, \rho f')',
\end{equation}
and Sturm-Liouville theory. Here,  
$D(L_w) \subseteq  \{ f \in \Ltwo \mid 
L_w f \in \Ltwo \}$ is a domain where $L_w$ is well-defined.
}
 
Indeed, the spectral problem \eqref{eq:PoincareWeakProblem} is formally equivalent to the eigenvalue problem
$	- L_w(f) = \eigval f $ with Neumann boundary conditions: $(\weight \, \rho f')(x) = 0$ if $x=a, b$. This comes from the following
integration by parts
\[\langle -L_wf,g\rangle=-\int_a^b (w\,\rho \,f')'\,g=\int_a^b (w\,\rho \,f')\,g'=\langle f',g'\rangle_w.\]
We refer to \cite{poincareintervals} for a rigorous proof in the case $\weight \equiv 1$ and $\rho > 0$ on a compact interval $[a, b]$. 
This eigenvalue problem can be rewritten in the Sturm-Liouville form, \begin{equation} \label{eq:SturmLiouvilleProblem}
	- (p f')' + q f = \eigval r f, \qquad p(a) f'(a) = p(b) f'(b) = 0
\end{equation}
with $p = \weight \, \rho$ and $q = 0$. If $a$ and/or $b$ are infinite, the boundary conditions above are interpreted as when taking the limit $a \to -\infty$ and/or $b\to +\infty$.\\
 
The Poincaré basis exists for a wide class of probability measures and weight functions. 
The next proposition gives sufficient conditions of existence.

\begin{prpstn}[Existence of the Poincaré basis]
\label{prop:sufficient_conditions_for_existence}
Let $\p = \weight \, \rho$. The Poincaré basis exists if at least one of the following conditions is verified:
\begin{itemize}
    \item [(i)] $1/\p$ belongs to $\LoneLeb$ 
    \item [(ii)] The anti-derivatives of $1/\p$ belong to $\Ltwo$. 
\end{itemize}
\end{prpstn}

Condition $(i)$ is more convenient than $(ii)$, but not always satisfied. For instance when $\mu$ is the uniform distribution on $(-1, 1)$ and $\weight(x) = 1 - x^2$, 
the corresponding Poincar\'e basis exists and is equal to the family of Legendre polynomials; in that case Condition $(ii)$ is verified \cite[Chapter 14, page 277]{Zettl} but Condition $(i)$ is not.
An example where the Poincaré basis does not exist is when $\mu$ is the uniform distribution on the interval $(-1,1)$ and $w(x)=(1-x^2)^2$ (see \cite{HerediaWeightPoincare}, \S 3.5, example of the Beta distribution with $\beta=1$).\\

The interest of the Poincaré basis for gradient-enhanced problems largely comes from the fact that its derivatives also form an orthogonal basis. This is actually a characteristic property of the Poincaré basis, as stated in the next proposition. \rev{We restrict our analysis to Condition (i) of Proposition~\ref{prop:sufficient_conditions_for_existence}, which suffices to guarantee the existence of the basis.} 

\begin{prpstn}[Poincaré basis and stability by differentiation]
\label{prop:PoincareBasis_and_Derivation}
	Let $\mu \in \probClass$, $\weight \in \funClassWeight$ and assume that $ 1/\p  \in \LoneLeb$ with $\p = \weight\,\rho$. Then, 
    \begin{enumerate}
        \item The Poincaré basis is an orthogonal basis of $\HoneW$ with $ \HoneNorm{\eigfun_\indBasis} = \sqrt{1+\eigval_\indBasis}$. Furthermore, the basis derivatives $(\eigfun'_\indBasis)_{\indBasis \geq 1}$ form an orthogonal basis of $\LtwoW$ with $ \LtwoNormW{\eigfun'_\indBasis}= \sqrt{\eigval_\indBasis}$.
        \item Conversely, if $(\psi_\indBasis)_{\indBasis \geq 0}$ is an orthonormal basis of $\Ltwo$ that belongs to $\HoneW$, with $\psi_0 \equiv 1$, and if $(\psi'_\indBasis)_{\indBasis \geq 1}$ is an orthogonal basis of $\LtwoW$, then $(\psi_\indBasis)$ is the Poincaré basis.
    \end{enumerate}
    \end{prpstn}

The proof is postponed to the appendix. It extends the one given in \cite{Poincare_Chaos_Sparse_2023} in which the simpler case $\weight \equiv 1$ is addressed. Note that in this reference, the condition $\psi_0 \equiv 1$, which is used in their proof and implicitly assumed in the context of chaos expansions, should be explicitly mentioned among the assumptions. It is mandatory since other bases could be obtained without it, e.g., by rotating the first two basis functions of the Poincaré basis and leaving the other ones unchanged.
Note also that if there exists other orthonormal bases (with $\psi_0 \equiv 1$) such that the basis function derivatives form an orthogonal system, 
then this system will not be complete, as a consequence of Proposition~\ref{prop:PoincareBasis_and_Derivation}, which is the case of the Fourier basis (see \cite{Poincare_Chaos_Sparse_2023}).\\ 

\begin{figure}[htbp]
\includegraphics[width=0.49\textwidth,height=5cm, trim=0 1cm 12cm 1cm,clip=TRUE]{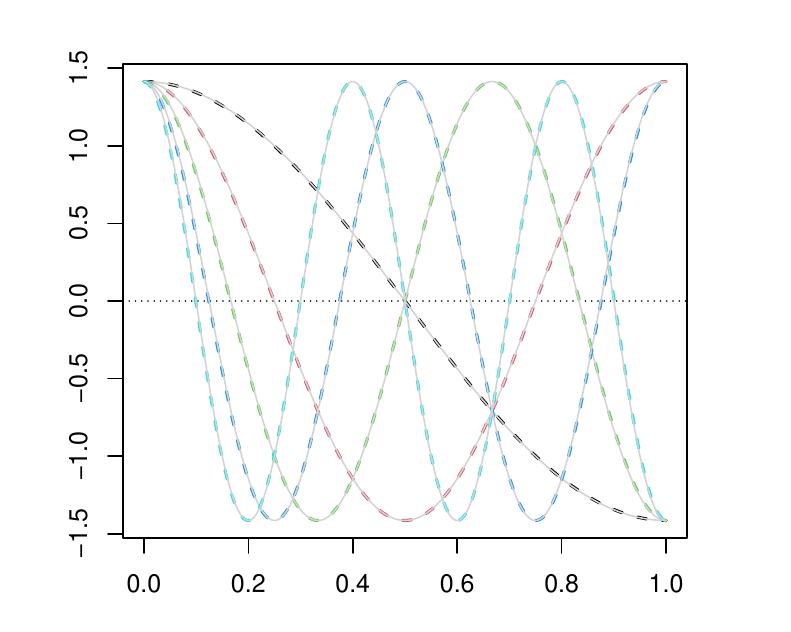}
\includegraphics[width=0.49\textwidth,height=5cm, trim=0 1cm 12cm 1cm,clip=TRUE]{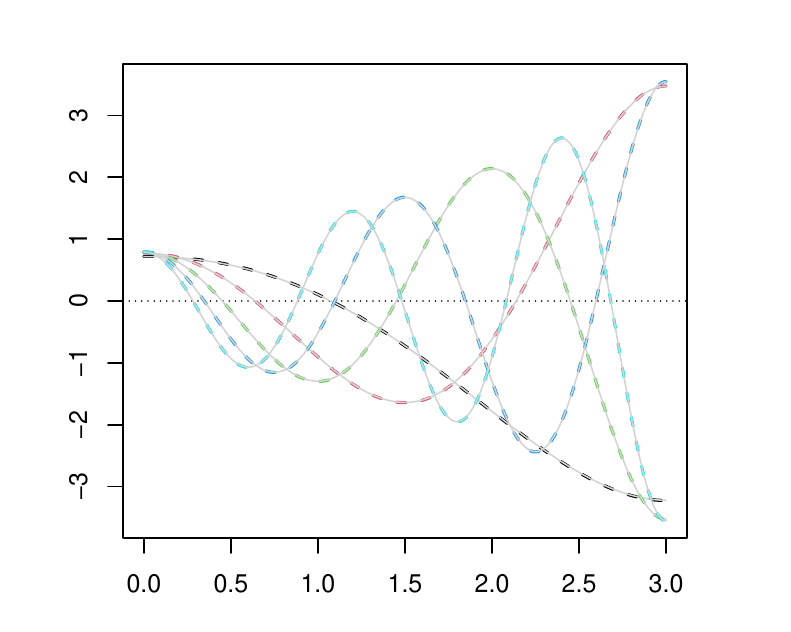}
\caption{First Poincaré basis functions (omitting the constant one) for $\mathcal{U}(0,1)$ and $\mathcal{E}(1)$ truncated on $[0,3]$, and $\weight \equiv 1$. Solid line: the basis function computed from the analytic expression; Dotted line: the basis function estimated by finite elements.}
\label{fig:PoincareBasis}
\end{figure}

The functions of the Poincaré basis verify oscillatory properties, well-known in Sturm-Liouville theory. Thus, the $\indBasis$-th basis function has exactly $\indBasis$ zeros \cite[Theorem 4.3.1, item (6)]{Zettl}.
As an illustration, we plot the first Poincaré basis functions for the uniform distribution on $[0,1]$ and the truncated exponential measure, truncated on $[0, 3]$, with a constant weight $\weight \equiv 1$ (Figure \ref{fig:PoincareBasis}).
Although such property is shared by orthogonal polynomials, Poincaré basis and orthogonal polynomials only coincide in the cases of the 
Normal, Gamma and Beta distributions, for a
unique choice of weight each time, which corresponds respectively to Hermite, Laguerre and Jacobi orthogonal polynomials (\cite{BGL}, \S 2.7). 
Apart from very specific choices of $\mu, \weight$ where the Poincaré basis is given explicitly, it must be computed numerically. This can be done efficiently by applying the finite element method to \eqref{eq:PoincareWeakProblem}, as detailed in \cite{poincareintervals, HerediaWeightPoincare}. This method has been used in Figure~\ref{fig:PoincareBasis} where we can see that the estimated basis functions are superimposed with the theoretical ones (here available). 

\rev{\subsection{Approximation properties for a bounded interval}}
\label{sec:approx_properties}

\rev{The next proposition provides approximation properties of the Poincaré basis, which are closely related to existing results for singular Sturm-Liouville problems \cite[Section 5.2.2]{canuto2006spectral}. As a first step, we restrict our analysis to a bounded interval $(a,b)$ and assume Condition (i) of Proposition~\ref{prop:sufficient_conditions_for_existence}, which suffices to guarantee the existence of the basis.

We will use the following notations. For a given $f \in \HoneW$, let $ \hat{f}_\indBasis := \langle f, \eigfun_\indBasis \rangle$ denote the coefficient of $f$ in the Poincaré basis,
and $P_N f = \sum_{j=0}^N \hat{f}_\indBasis \eigfun_\indBasis $ 
the approximation of $f$ in the finite-dimensional space spanned by the first $N$ basis functions. Finally, let $\Hm{m}$ denote the usual Sobolev space of order $m$, with norm defined by
$\Vert f \Vert_{\Hm{m}}^2 :=
\sum_{i=0}^m \Vert  f^{(i)} \Vert^2$.\\
}

\rev{
\begin{prpstn}[Approximation property of the Poincaré basis]
\label{prop:approximation_property}
Let $(a,b)$ be a bounded interval.
Let $\mu \in \probClass$ with density $\rho$ and let
$\weight \in \funClassWeight$. Denote $p=wr$.
Let $k \in \N^*$  and 
assume:
\begin{itemize}
    \item \emph{Smoothness condition:} the functions $w$ and $\tau:=(w\,\rho)'/\rho$ are of class $\Cm{2k-2}$ on $[a,b]$.
    \item \emph{Basis existence condition:} $1/p \in \LoneLeb$
    \item \emph{Singular Sturm-Liouville problem:} $p(a)=p(b)=0$
\end{itemize}
Then the following assertions hold.
\begin{enumerate}
    \item[(i)] The iterated operator $(-L_w)^k$ is continuous from $\Hm{2k}$ to $\Ltwo$. Moreover, for  all $f\in \Hm{2k}$ and all $j \geq 1$ we have the identity
    \begin{equation} \label{eq:coefFormulaWithLiterate}
 \hat{f}_\indBasis 
= \frac{1}{\eigval^k_\indBasis} \langle (-L_w)^k f, \eigfun_\indBasis \rangle.
\end{equation} 
    \item[(ii)] There exists a constant $C>0$ such that for all $f\in \Hm{2k}$ and all $j \geq 1$ 
\begin{equation} \label{eq:coefDecrease}
    \vert \hat{f}_\indBasis \vert \leq \frac{C}{\indBasis^{2k}} \Vert (-L_w)^{k} f\Vert.
\end{equation}
    \item[(iii)] There exists a constant $C_k>0$ such that for all $f\in \Hm{2k}$ and all $N \geq 1$  
    \begin{equation} \label{eq:approxError}
        \Vert  f - P_N f \Vert \leq \frac{C_k}{N^{2k}} \Vert  f \Vert_{\Hm{2k}}.  
    \end{equation}
\end{enumerate}
\end{prpstn}
}

\rev{
Proposition~\ref{prop:approximation_property} shows that, when $p$ vanishes at the boundaries, then, under technical conditions related to the smoothness of $w, r$, the ($L^2$) approximation error $\Vert f - P_N f \Vert $ decreases at the rate $N^{-m}$ for all functions $f$ in $\Hm{m}$ and all even integer $m$, which is the rate observed for the best polynomial approximation (for the $L^\infty$ norm) of a function of class $\Cm{m}$ (see e.g., \cite[Chapter 7]{trefethen2013approximation} or \cite[Chapter 4, \S 6, `Jackson's Theorem V', (iii)]{cheney1966approximation}). 
The condition that $p$ vanishes at the boundaries is essential and corresponds to singular Sturm Liouville problems. If it is not verified, then Equation~\ref{eq:coefFormulaWithLiterate}, which is obtained by integration by part, must be modified to involve the derivatives of $f$ at the boundaries. The speed of convergence is then slower or limited to functions of $\Hm{m}$ for which all derivatives up to order $m$ vanish at the boundaries. We refer to \cite{xiu2010numerical} and \cite{canuto2006spectral}, Section 5.2., 
for a more detailed discussion of this approach
and its application to Legendre polynomials, which constitute a particular instance of a Poincar\'e basis. 
We also refer to 
\cite{adcock2009univariate},
in which the case where $\mu$ is uniform, with constant weights, is investigated (`modified Fourier series').

This result gives some guidance to choose the weight $w$. Indeed, if the pdf $r$ of $\mu$ does not vanish at the boundaries, then it may be useful to choose $w$ that vanishes at the boundaries so that $p = w r$ also vanishes at the boundaries.}

\subsection{A special choice of weight
\texorpdfstring{{\boldmath$\wlin$}}{wlin}}
\label{sec:weight_w_lin}
A desirable property in practice is that the Poincar\'e basis includes linear functions, as they often provide good approximations on the behavior observed in practical models.
Since the basis is completely determined by the pair of probability measure $\mu$ and weight $\weight$, enforcing the second eigenfunction to be linear is achieved by appropriately selecting the weight. 
\rev{Assuming that $\int_a^b |x|\,d\mu(x)<\infty$} and denoting the mean of $\mu$ by $m=\int_a^b \rev{y \,d\mu(y)}$, such a weight choice is explicitly defined as
\begin{equation}
\label{eq:wlin}
    \wlin(x)=-\frac{1}{\rho(x)}\int_a^x (y-m)\, \rho(y)\,dy,\quad \mbox{for all }x\in (a,b).
\end{equation}
\rev{
This weight coincides with the so-called Stein kernel of the measure $\mu$, which, since $\mu$ is assumed to have a finite first moment, is uniquely defined in the present one-dimensional setting (see \textit{e.g.} \cite{saumard}). Moreover, this weight} is exactly the one associated with the only three cases where the Poincar\'e basis is polynomial (Hermite, Laguerre and Jacobi).
As such, for these particular cases, it has a classical closed-form expression and has already been used for gradient-enhanced surrogate modeling (see e.g. \cite{Adcock_2019, Peng2016, Guo2018}).
\rev{Notice that when $(a,b)$ is a bounded interval, $(\wlin \, r)$ vanishes at the boundary. Thus, under additional conditions on $r$, the Poincaré basis has the strong approximation properties of Proposition~\ref{prop:approximation_property}. This happens for instance when $r$ does not vanish on $[a,b]$ and is of class $C^{2k-2}$, provided that $1/p \in L^1(a,b)$.} 

\begin{figure}[htbp]
\includegraphics[width=0.49\textwidth]{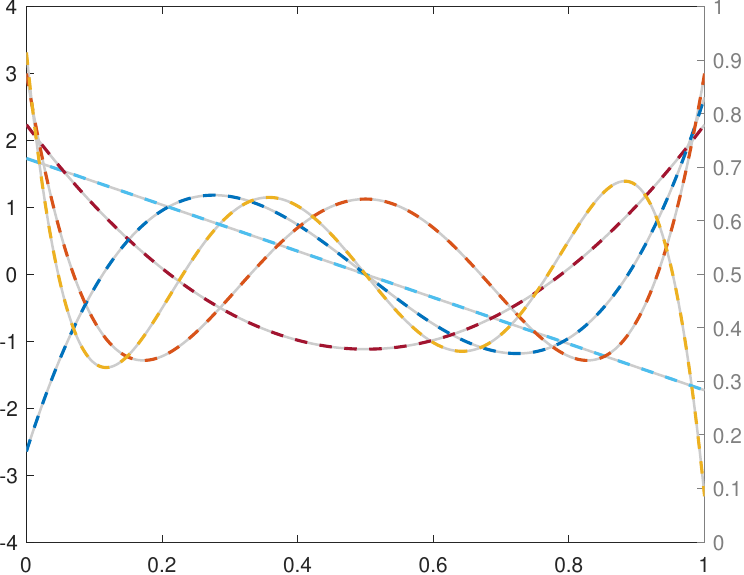}
\includegraphics[width=0.49\textwidth]{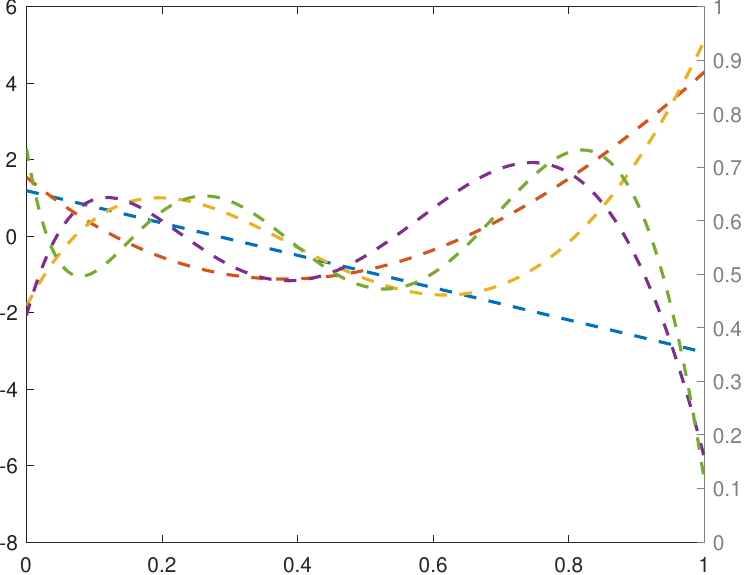}
\caption{First Poincaré basis functions (omitting the constant one) for $\mathcal{U}(0,1)$ and $\mathcal{E}(1)$ truncated on $[0,3]$, and $\weight = \wlin$. The dotted lines represent the basis functions estimated by finite elements. On the left plot, the solid lines are the theoretical functions associated to $\mathcal{U}(0,1)$ (\textit{i.e.} the Legendre polynomials)}
\label{fig:PoincareBasis_wlin}
\end{figure}

In the GSA context, the weight $\wlin$ was first considered in \cite{Song}, where the authors develop some applications for models involving probability measures for which $\wlin$ is computed explicitly. It was further investigated in \cite{HerediaWeightPoincare}, where it is shown that this weight is particularly well suited for models exhibiting linear trends. \cite{HerediaWeightPoincare} also provides a numerical method to approximate $\wlin$, enabling its practical use for any probability measure $\mu\in \probClass$. The idea of such a method is simple. One first solves numerically the  Cauchy problem that we obtain after multiplying $\rho$ with both sides of  \eqref{eq:wlin}, followed by differentiation:
\begin{equation}
\label{eq:cauchy_problem_wlin}
       \left\{\begin{array}{rl}
        (\wlin\, \rho)'(x)&=-(x-m)\,\rho(x) \quad \mbox{on }(a,b),\\
        (\wlin\, \rho)(a)&=0,
       \end{array}\right.
\end{equation}
and then divide the approximated solution of $\wlin\, \rho$ by $\rho$. This approach avoids using numerical integration to compute the anti-derivative in \eqref{eq:wlin} for multiple values $x\in (a,b)$,
which would imply a higher computational cost.
The resulting approximated weight is very accurate as seen in \cite{HerediaWeightPoincare}, at least when using the Runge-Kutta~4 method to solve \eqref{eq:cauchy_problem_wlin} (see \textit{e.g.} \cite{BurdenNumerical}). Then the Poincar\'e basis functions are also well approximated, using a finite element discretization as mentioned in the preceding section. This is illustrated in the left plot in Figure~\ref{fig:PoincareBasis_wlin} in the case of the uniform distribution, whose eigenfunctions coincide with the Legendre polynomials. The plot on the right displays the approximated eigenfunctions associated to a truncated exponential distribution, for which closed-form expressions are not available.

\section{Poincaré chaos expansions}
\label{sec:PoinCE}

In this section, we consider the chaos expansion obtained from the (univariate) Poincaré bases, called Poincaré chaos expansion. Thus, for all $\indVar=1, \dots, d$, we denote by $(\psi_{\indVar, \indBasis})_{\indBasis \in \N}$ the Poincaré basis associated to $\mu_\indVar$, and $(\eigval_{\indVar, \indBasis})_{\indBasis \in \N}$ the corresponding eigenvalues. The chaos expansion is then defined by using Eq.\eqref{eq:CEdef}.

\subsection{DGSM computation from Poincar\'e chaos expansion}
\label{sec:DGSMfromPoincare}

The orthogonality property of the Poincaré basis derivatives can be leveraged to obtain a closed-form expression of DGSM.

\begin{prpstn}
\label{prop:DGSMfromPoinCE}
Let $\indVar \in \{1, \dots, d \}$. Then the DGSM associated to $X_\indVar$ can be simply derived from the Poincaré chaos expansion as
\begin{equation} \label{eq:DGSMfromPoinCE}
\nu_\indVar = \sum_{\alp \in \N^d, \alpha_\indVar \geq 1} \eigval_{\indVar, \alpha_\indVar} c_\alp^2.
\end{equation}
Furthermore, $S_\indVar^\tot \leq C_P(\mu_\indVar, \weight_\indVar) \frac{\nu_\indVar}{\Var(\cm(X))}$, where $C_P(\mu_\indVar, \weight_\indVar) = 1/\eigval_{\indVar, 1}$ is the Poincaré constant of $\mu_\indVar$ for the weight $\weight_\indVar$.
\end{prpstn}

The originality of the proposition concerns the first result, which has been proved so far only in the case $\weight_\indVar \equiv 1$ \cite{Poincare_Chaos_Sparse_2023}. The second part, linking the Sobol' indices to DGSM, actually does not require the existence of the Poincaré basis, and was proved in \cite{HerediaWeightPoincare} for a general weight, by using a different approach, following earlier results from \cite{sobkuc09, Lamboni_et_al_2013} in the case $\weight_\indVar \equiv 1$.

\rev{

\subsection{Sparse regression}
\label{sec:sparse_regression}

To compute the coefficients of a chaos expansion, it needs to be truncated to a finite set of basis elements.
Based on a finite number of model evaluations, the coefficients can then efficiently be computed by sparse regression, a least squares regression technique which enforces sparsity in the chaos coefficients through regularization \cite{LuethenSIAMJUQ2020}. 

When model derivatives are available, it is desirable to include this information in the regression process. In the following, we describe two ways to achieve this. More details can be found in Appendix~\ref{sec:app:regression}.

Note that the property of the basis having orthogonal partial derivatives is crucial for these methods. Although one may attempt to perform gradient-informed regression using chaos expansions without orthogonal partial derivatives (such as unweighted or non-classical families of polynomials), this approach is not straightforward, as pointed out in \cite{Guo2018}. The authors even illustrate for Legendre PCE how directly including derivative information (without weighting/preconditioning) can deteriorate the stability of the regression matrix.

}

\rev{
\subsubsection{Aggregating derivative expansions}
\label{sec:averaged_deriv_expansions}

One possible approach, presented in \cite{Poincare_Chaos_Sparse_2023} and detailed in Appendix~\ref{app:averaged_deriv_expansions}, is to use each of the $d$ sets of partial derivatives separately to compute a subset of coefficients by sparse regression. 
This subset is smaller than the original set of coefficients, because differentiation with respect to a variable removes all basis elements which are constant in this variable.
As a post-processing step, coefficients computed by more than one derivative expansion are finally averaged to yield a single estimate. 

}

\rev{
\subsubsection{Gradient-enhanced regression (combined regression)}
\label{sec:gradient_enhanced_regression}

Model evaluations and gradient information can also be collected into a single tall regression problem \cite{Adcock_2019}, see Appendix~\ref{app:gradient_enhanced_regression} for details.
Here the basis function and basis partial derivative evaluations are stacked, resulting in a system with as many columns as there are basis functions, and $N(d+1)$ rows, with $N$ the number of experimental design points and $d$ the number of input variables.
It is crucial to use weighted or unweighted Poincar\'e basis functions, and normalize the columns of the matrix, to ensure coherence and isotropy which are needed for accuracy and stability \cite{Hampton2015}.

}

\section{Numerical experiments}
\label{sec:numExp}

We now present some numerical experiments demonstrating the performance of gradient-enhanced Poincar\'e chaos expansions for sensitivity analysis and surrogate modeling.

\subsection{Implementation}
We have implemented the methodology within the UQLab framework \cite{MarelliUQLab2014}. 
We consider two different settings for the basis functions:
\begin{itemize}
    \item Fully unweighted case: For each variable $X_\indVar$, we use the one-dimensional Poincar\'e basis associated to the constant weight $\weight_\indVar\equiv 1$.
    \item Fully weighted case: for each variable $X_\indVar$, we set $\weight_\indVar$ equal to $\wlin$ introduced in Section~\ref{sec:weight_w_lin}, and use the associated Poincar\'e basis. The procedure to compute $\wlin$  is detailed in that section.
\end{itemize}

\noindent
For simplicity, we denote by $\ve{w}$ the collection of all the one-dimensional weights $\weight_\indVar$. The Poincar\'e bases are approximated using a finite element discretization of the spectral problem \eqref{eq:PoincareWeakProblem}, as presented in \cite{poincareintervals}. 

We consider the following chaos expansions:
\begin{itemize}
    \item The standard expansion based exclusively on basis functions (without derivatives), given in \eqref{eq:chaosExpansion}. We abbreviate it as PoinCE and wPoinCE in the unweighted and weighted cases, respectively. 
    \rev{In certain cases, these expansions coincide with PCE (see Section \ref{sec:weight_w_lin}).}
    \item The averaged derivative expansions introduced in \rev{Section~}\ref{sec:averaged_deriv_expansions}, here abbreviated as PoinCE-der-aggr and wPoinCE-der-aggr. 
    \item The expansions obtained from the gradient-enhanced regression in Section \ref{sec:gradient_enhanced_regression}. 
    We abbreviate them as PoinCE-comb-regr and wPoinCE-comb-regr.
\end{itemize}
The sparse regression problems are solved using the LARS solver available in UQLab \cite{UQdoc_21_104}, which achieves model selection based on the leave-one-out error of the associated surrogate model. We use a total degree truncation set of the specified degree.
\rev{The experimental design is a Latin Hypercube design sampled from the joint pdf of the input random variables.}

We compare the performance of each expansion with respect to the following quantities:
\begin{itemize}
    \item The $\HoneW$ error on a validation set  given by
    \begin{align}
E_{\xHone} &= \hat{\mathbb{E}}[(\cm(\cx_\text{val}) - \cm^\text{surr}_\ca(\cx_\text{val}, \ve c))^2]  
\nonumber\\
\label{eq:HoneW}
& + \sum_{\indVar = 1}^d \hat{\mathbb{E}}\left[
w_\indVar\bigr(\cx_{\text{val},\indVar}\bigr)
\left(\partderiv[\cm]{\indVar}(\cx_\text{val}) - \partderiv[\cm]{\indVar}\cm^\text{surr}_\ca(\cx_\text{val}, \ve c)\right)^2\right].
\end{align}
\item \rev{The $\Ltwo$ error on a validation set, defined by the first right-hand side term of Eq.~\eqref{eq:HoneW}.}
 \item 
 The total Sobol' indices, computed from Eq.\eqref{eq:totalSobolFromCE}. Since the exact expressions of the indices are unavailable for the models we consider, for the sake of comparison we include high-precision estimates obtained via PCE with a large sample size, that we refer to as the ``true'' values.
\end{itemize}

    The numerical experiments are performed multiple times using different input sample sizes, that we refer to as experimental design sizes (ED sizes). For each experiment, both the model and its derivatives are evaluated at the same input points. In addition, we perform $30$ bootstrap replicates of each estimation of total Sobol' indices, \rev{$\Ltwo$} and $\HoneW$ errors. They are displayed with boxplots to represent confidence intervals.

\subsection{Toy model with interaction}
\label{sec:Toy_model_with interaction}
Consider the model
\begin{equation}
    f(\ve X) = \prod_{\indVar=1}^d \frac{d/4}{d/4 + (X_\indVar - a_\indVar)^2},\quad \quad a_\indVar = \frac{(-1)^\indVar}{\indVar + 1}, 
\end{equation}
where $d$ is the dimension and $\ve X=(X_1,\dots,X_d)$ is a vector of independent uniform random variables $X_\indVar \sim \cu(-1,1)$. As such, the Poincar\'e basis $(\psi_\indBasis)$ is known in both the unweighted and the weighted case. When $\weight_\indVar\equiv 1$ the basis functions are the orthonormal cosines $\psi_\indBasis(x)=\sqrt{2}\cos(\indBasis\pi x)$ and when $\weight_\indVar=\wlin$ they coincide with the Legendre polynomials. Hence wPoinCE is actually PCE.

This model was used by \cite{Adcock_2019} for gradient-enhanced regression with Legendre polynomials and points drawn from the uniform density, as we do in the weighted case, as well as with Chebyshev polynomials and points drawn from the Chebyshev density. 
They look at the three different cases $(d,s) \in \{(4,72),(8,23),(12,14)\}$, where $s$ is the degree of the hyperbolic cross index set they use for basis truncation. Here we consider $d = 4$ dimensions and a total degree truncation basis of $p = 8$.  Besides these precisions, the main difference in the implementation between their approach and ours (in the weighted case) is that they use the SPGL1 solver \cite{SPGL1} instead of LARS.

\begin{figure}[htbp]
  \centering
  \subfloat[$\HoneW$ error\label{fig:AdcockModel:PoinCE:H1}]{
    \includegraphics[width=0.31\textwidth]{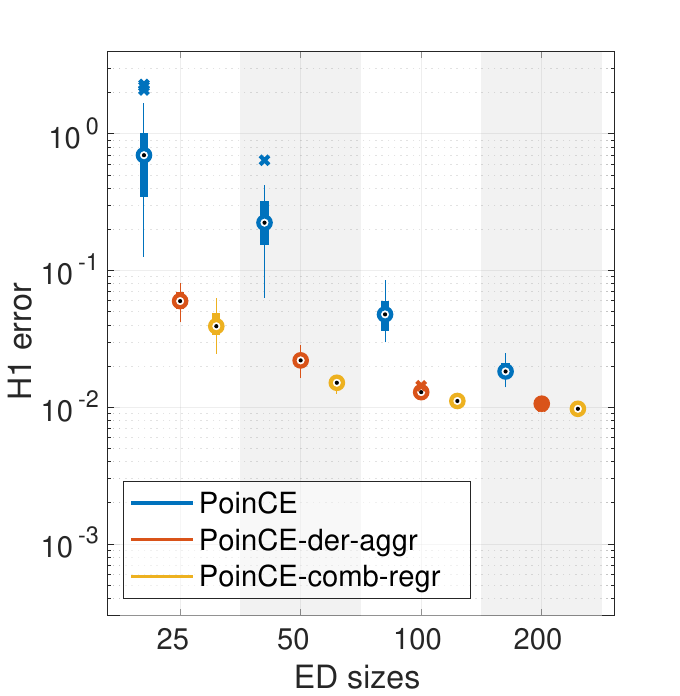}
  }\hfill
  \subfloat[Total Sobol' index of $X_1$\label{fig:AdcockModel:PoinCE:SobolX1}]{
    \includegraphics[width=0.31\textwidth]{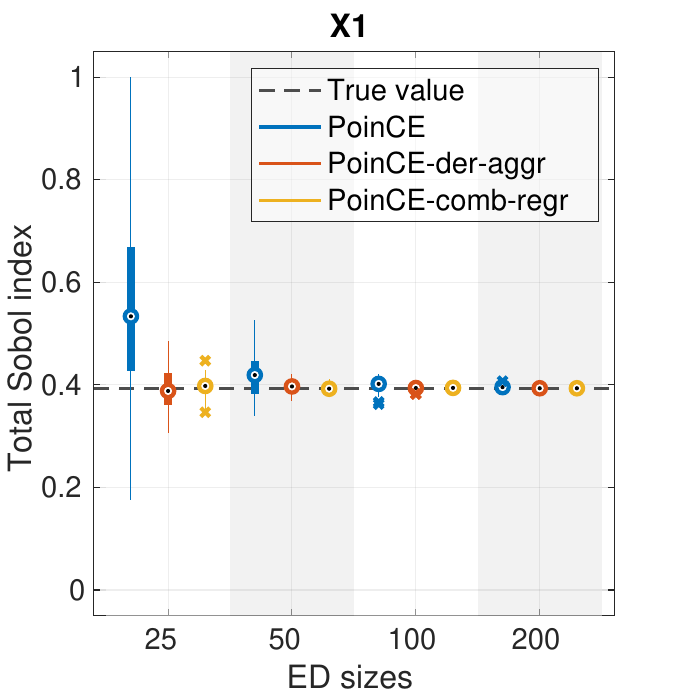}
  }\hfill
  \subfloat[Total Sobol' index of $X_4$\label{fig:AdcockModel:PoinCE:SobolX4}]{
    \includegraphics[width=0.31\textwidth]{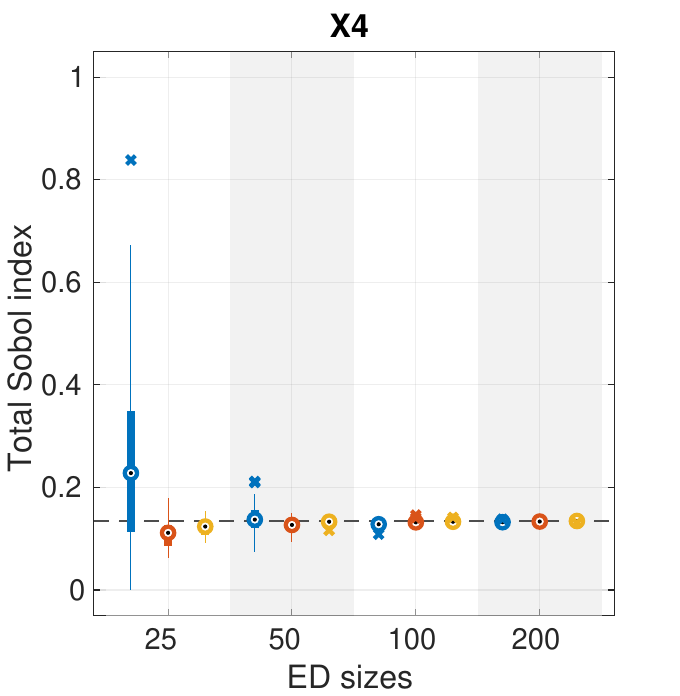}
  }
  \caption{Results for the toy model in the unweighted case.
    Dotted lines represent the true total Sobol' indices.}
  \label{fig:AdcockModel:PoinCE}
\end{figure}

\begin{figure}[htbp]
  \centering
  \subfloat[$\HoneW$ error\label{fig:AdcockModel:wPoinCE:H1}]{
    \includegraphics[width=0.31\textwidth]{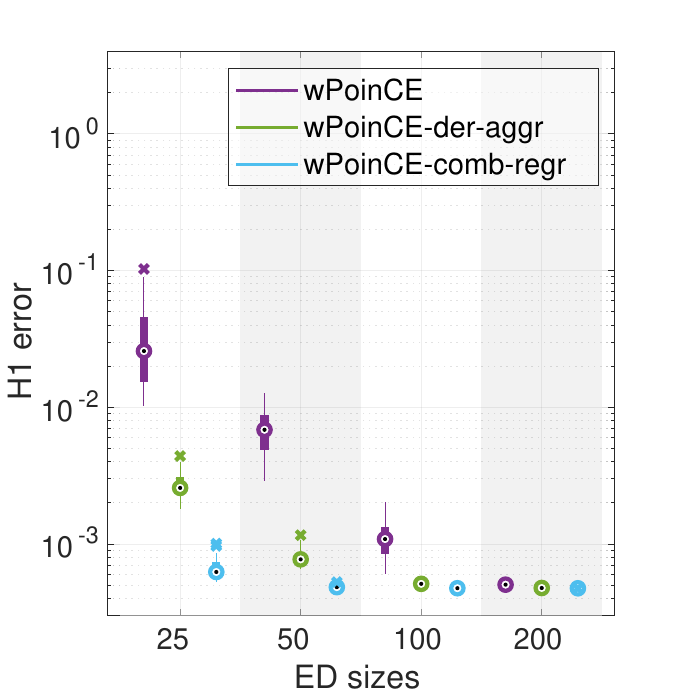}
  }\hfill
  \subfloat[Total Sobol' index of $X_1$\label{fig:AdcockModel:wPoinCE:SobolX1}]{
    \includegraphics[width=0.31\textwidth]{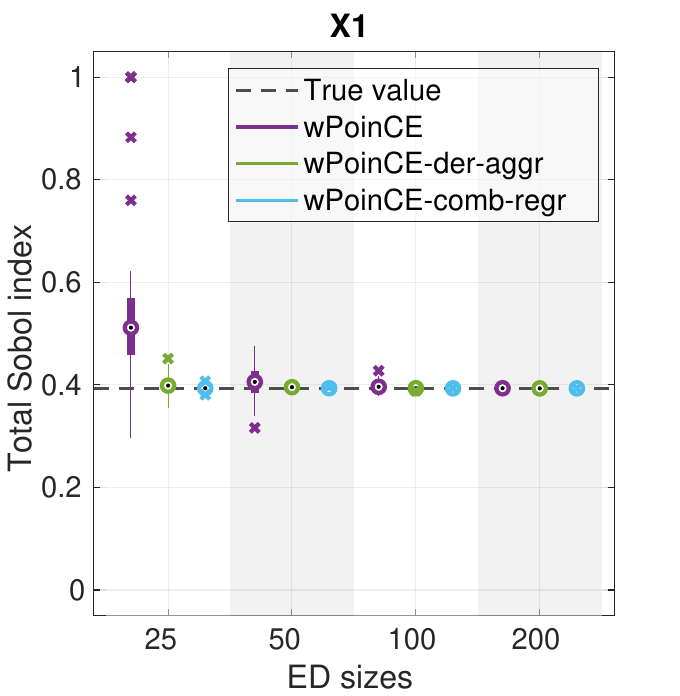}
  }\hfill
  \subfloat[Total Sobol' index of $X_4$\label{fig:AdcockModel:wPoinCE:SobolX4}]{
    \includegraphics[width=0.31\textwidth]{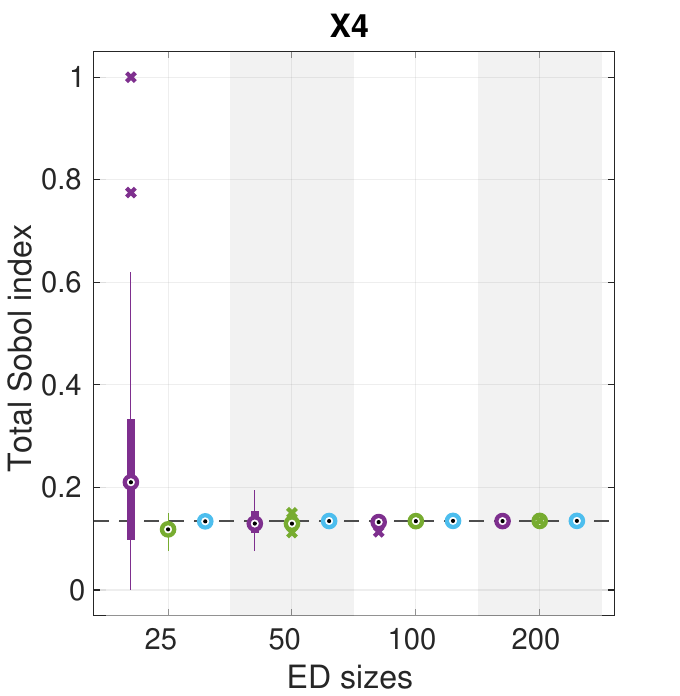}
  }
  \caption{Results for the toy model in the weighted case.
    Dotted lines represent the true total Sobol' indices.}
  \label{fig:AdcockModel:wPoinCE}
\end{figure}

Figure \ref{fig:AdcockModel:PoinCE} presents the results in the unweighted case. It displays estimations of both the $\HoneW$ errors and the total Sobol' indices of the variables $X_1$ and $X_4$. Observe that for both errors and indices, the expansions based on derivative evaluations (PoinCE-der-aggr and PoinCE-comb-regr) outperform PoinCE, which only relies on model evaluations. Moreover, the two derivative-based expansions are comparable in this model, with a slight advantage in favor of PoinCE-comb-regr. In the weighted case, the same situation occurs with the corresponding expansions, as shown in Figure \ref{fig:AdcockModel:wPoinCE}. Actually, the advantage of using wPoinCE-comb-regr is more evident here, in terms of the $\HoneW$ error.

The conclusions above make sense if we suppose that all the model and gradient evaluations are already available, \emph{i.e.} that the derivatives involve no additional cost. To compare the methods more realistically, one may define an equivalence rate between model and gradient evaluations. For instance in \cite{Adcock_2019,Jakeman2015} it is supposed that evaluating the gradient once is equally expensive as evaluating the model.
When adapting this equivalence rule, our conclusions remain valid (for example in the unweighted case, each blue boxplot is compared to the red and yellow ones on the left next to it). Furthermore, our weighted results involving Legendre polynomials, show improved $\HoneW$ errors compared to those reported in \cite{Adcock_2019}. This suggests an advantage of using LARS over the SPGL1 solver, and confirming the findings of \cite{LuethenSIAMJUQ2020} on other toy functions. 

Another natural equivalence rule to consider is the following: the cost of one gradient evaluation equals the cost of four model evaluations (since the gradient has $d=4$ partial derivatives). Although Figures \ref{fig:AdcockModel:PoinCE} and \ref{fig:AdcockModel:wPoinCE} do not illustrate this exact cost-equivalence setting, since some ED sizes should be $1+d=5$ larger than others, they provide a very similar scenario. By comparing the results using ED sizes of $25$ and $50$ with those of $100$ and $200$, respectively, we conclude that PoinCE (resp. wPoinCE) becomes comparable to PoinCE-comb-regr (resp. wPoinCE-comb-regr).

\begin{figure}[htbp]
    \centering  \includegraphics[width=0.35\linewidth]{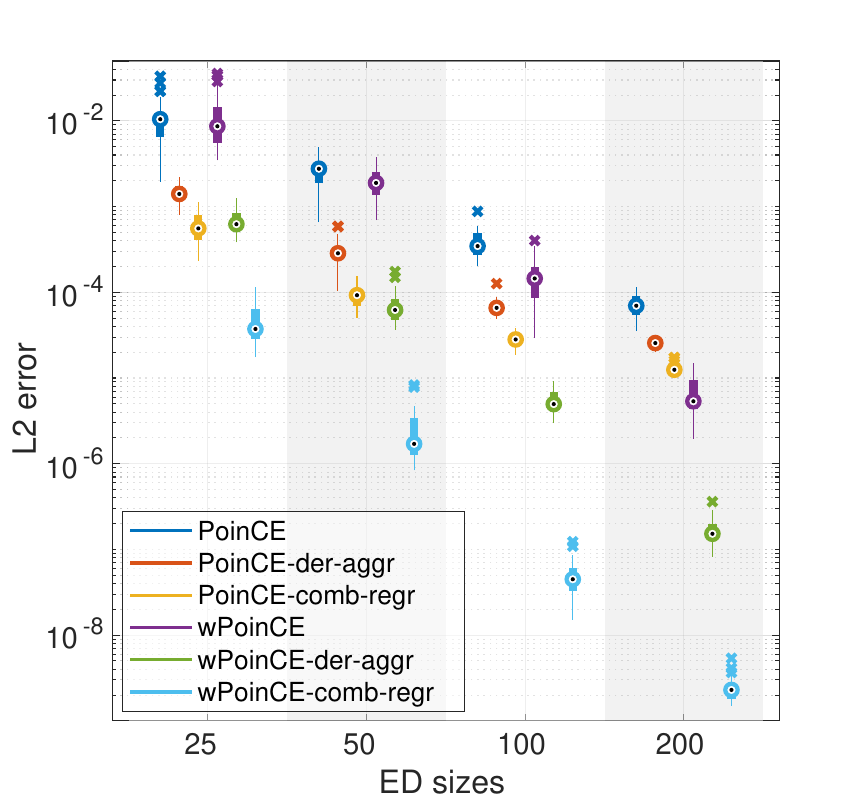}
    \caption{$L_2(\mu)$ error. Results for the toy model, both unweighted and weighted cases.    }
    \label{fig:AdcockModel:PoinCE:L2}
\end{figure}

We also compare the unweighted and weighted cases. First, in Figures \ref{fig:AdcockModel:PoinCE} and \ref{fig:AdcockModel:wPoinCE} we find similar approximations on the total Sobol' indices.
Second, we cannot draw conclusions about the surrogate models by comparing their $\HoneW$ errors, as these depend on $\weight$. 
We rather use the $\Ltwo$ errors. 
They are displayed in Figure \ref{fig:AdcockModel:PoinCE:L2}, showing a clear advantage in the weighted case. 

\subsection{Flood model}

We now consider a flood model commonly used to test GSA methodologies (see \textit{e.g.}
\cite{HerediaWeightPoincare,Poincare_Chaos_Sparse_2023,poincareintervals,PoincareChaos}). In this model, the maximal annual overflow of a river (measured in meters) is given by
\[
		S=Z_v-H_d-C_b+\left(\frac{Q}{B K_s}\sqrt{\frac{L}{Z_m-Z_v}}\right)^{\frac{3}{5}},
\]
and our primary interest lies in the annual maintenance cost of a dyke constructed next to it:
\[
C =\mathbb{1}_{S>0}+\left(0.2+0.8\left(1-e^{-\frac{1000}{S^4}}\right)\right)\mathbb{1}_{S\leq 0}+\frac{1}{20}\max\left\{ H_d,8\right \}.
\]

The inputs are supposed to be independent random variables, with details provided in Table \ref{tab:variables_flood}. Note that none of the variables, except $H_d$, follow the standard distributions: normal, exponential or uniform. Thus, for these variables, the one-dimensional Poincar\'e basis functions are not polynomials. We have justified the existence of their Poincar\'e bases in the unweighted case.
Indeed, in this setting, the Gumbel and normal truncated distributions satisfy condition $(i)$ in Proposition \ref{prop:sufficient_conditions_for_existence}, and any triangular measure satisf\rev{ies} condition $(ii)$.
However in the weighted setting with $\wlin$, we do not have closed-form expressions for the functions involved in conditions $(i)$ and $(ii)$.
Establishing the existence of the corresponding Poincar\'e bases for these measures requires a more detailed analysis, which is beyond the scope of this article. 


\begin{table}[h]
\begin{center}
\bgroup
\def\arraystretch{1.2}
\begin{tabular}{llll}
\hline
Input & Meaning & Unit & Probability measure   \\ \hline 
$Q$ & Max. flow rate & $m^3/s$ & Gumbel $\mathcal{G}(1013, 558)|_{[500, 3000]}$ \\  
$K_s$ & Strickler coefficient  & --- & Gaussian $\mathcal{N}(30,64)|_{[15,75]}$ \\ 
$Z_v$ & Downstream level & $m$ & Triangular $\mathcal{T}(49, 50, 51)$\\ 
$Z_m$ & Upstream level & $m$ & Triangular $\mathcal{T}(54, 55, 56)$ \\
$H_d$ & Dyke height & $m$ & Uniform $\mathcal{U}(7,9)$\\ 
$C_b$ & Bank height & $m$ & Triangular $\mathcal{T}(55, 55.5, 56)$ \\
$L$ & River length & $m$ &  Triangular $\mathcal{T}(4990, 5000, 5010)$ \\ 
$B$ & River width & $m$ & Triangular $\mathcal{T}(295, 300, 305)$ \\ 
\hline
\end{tabular}
\egroup
\caption{Input variables for the flood model. The notations $\mathcal{G}(\eta,\beta)$ ($\eta\in \R$, $\beta>0$) $\mathcal{T}(a,c,b)$ ($a<c<b$) are reserved to Gumbel and triangular distributions, respectively. The notation $|_I$ means that the distribution is truncated on the interval $I$.}
\label{tab:variables_flood}
\end{center}
\end{table}

\begin{figure}[h!]
  \centering
  \subfloat[$\HoneW$ error\label{fig:FloodModel:PoinCE:H1}]{
    \includegraphics[width=0.31\textwidth]{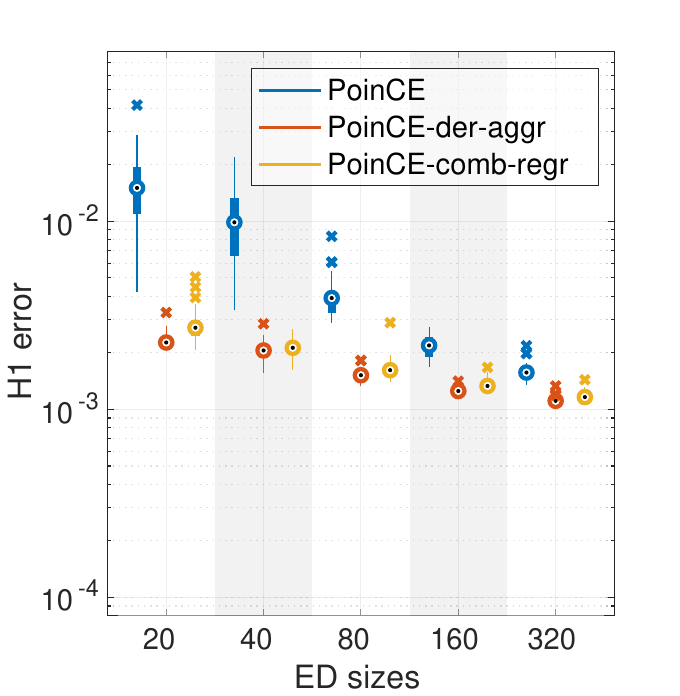}
  }\hfill
  \subfloat[Total Sobol' index of $Q$\label{fig:FloodModel:PoinCE:SobolX1}]{
    \includegraphics[width=0.31\textwidth]{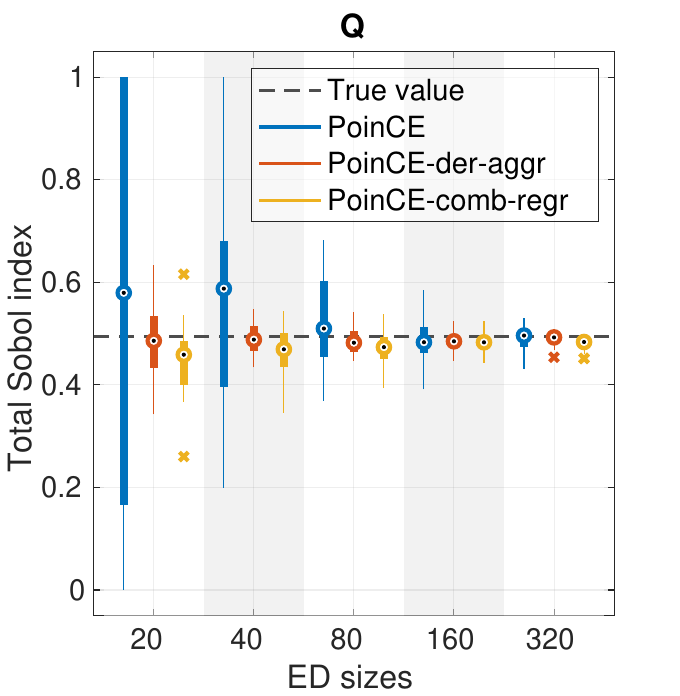}
  }\hfill
  \subfloat[Total Sobol' index of $C_b$\label{fig:FloodModel:PoinCE:SobolX6}]{
    \includegraphics[width=0.31\textwidth]{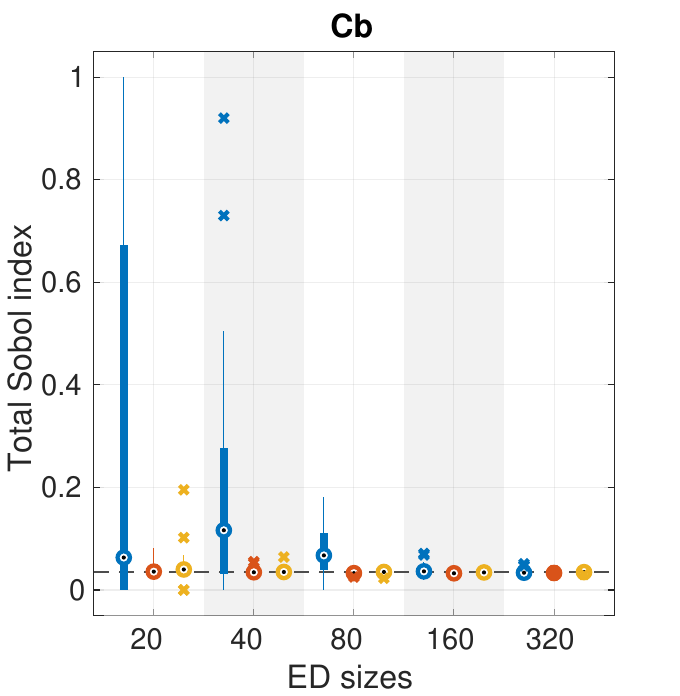}
  }
  \caption{Results for the flood cost model in the unweighted case.}
  \label{fig:FloodModel:PoinCE}
\end{figure}

Our results for the unweighted case are shown in Figure \ref{fig:FloodModel:PoinCE}, which displays the $\HoneW$ errors and total Sobol' indices estimations for one influential variable ($H$) and one non-influential variable ($C_b$). Here again, the methods based on model derivatives outperform PoinCE and in this model, PoinCE-der-aggr performs slightly better than PoinCE-comb-regr. In the weighted case we obtain similar conclusions, except that the best results are provided by wPoinCE-comb-regr. 

Since some of the methods leverage model derivatives while others do not, a fair comparison requires defining the cost associated with gradient evaluations.
If we follow the rule adopted in \cite{Adcock_2019,Jakeman2015}, where one gradient evaluation is assumed to have the same cost as one model evaluation, we observe that the methods relying on model derivatives perform better than those based on model evaluations only (for instance, in Figure \ref{fig:FloodModel:PoinCE}, each blue boxplot is compared to the red and yellow ones to the left next to it).
If, however, gradients are approximated using finite differences, then one
gradient evaluation costs $d=8$ times a single model evaluation, and the closest comparison is obtained by matching the results using ED sizes of 20 and 40 with those of 160 and 320, respectively.
Under this cost model, the performance of PoinCE and wPoinCE is comparable to that of the derivative-based expansions.

\begin{figure}[h!]
  \centering
  \subfloat[$\HoneW$ error\label{fig:FloodModel:wPoinCE:H1}]{
    \includegraphics[width=0.31\textwidth]{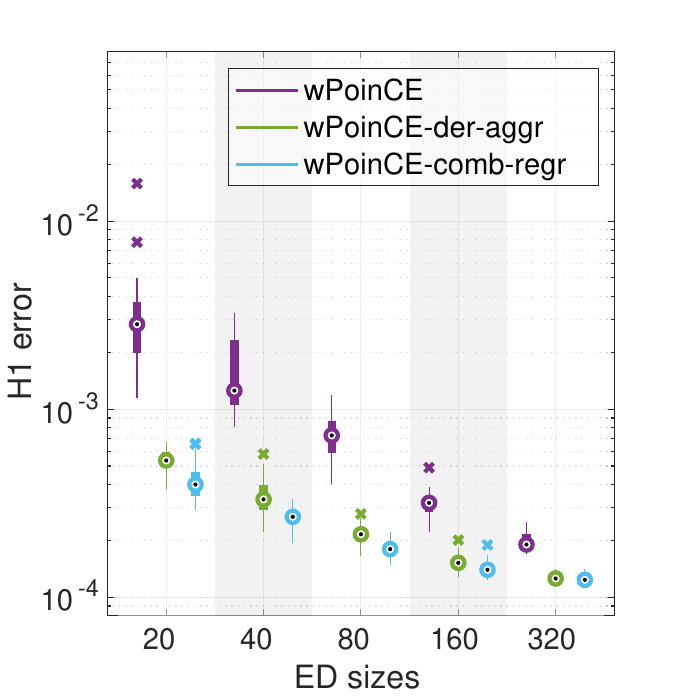}
  }\hfill
  \subfloat[Total Sobol' index of $Q$\label{fig:FloodModel:wPoinCE:SobolX1}]{
    \includegraphics[width=0.31\textwidth]{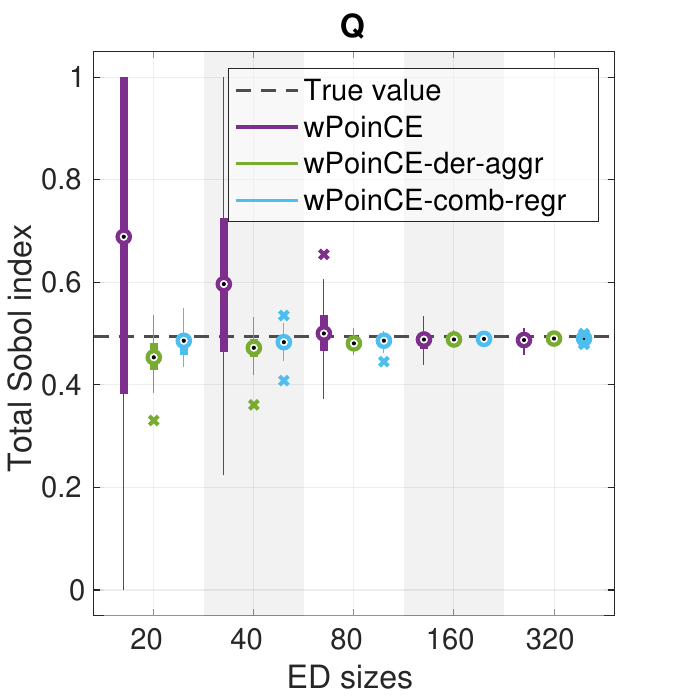}
  }\hfill
  \subfloat[Total Sobol' index of $C_b$\label{fig:FloodModel:wPoinCE:SobolX6}]{
    \includegraphics[width=0.31\textwidth]{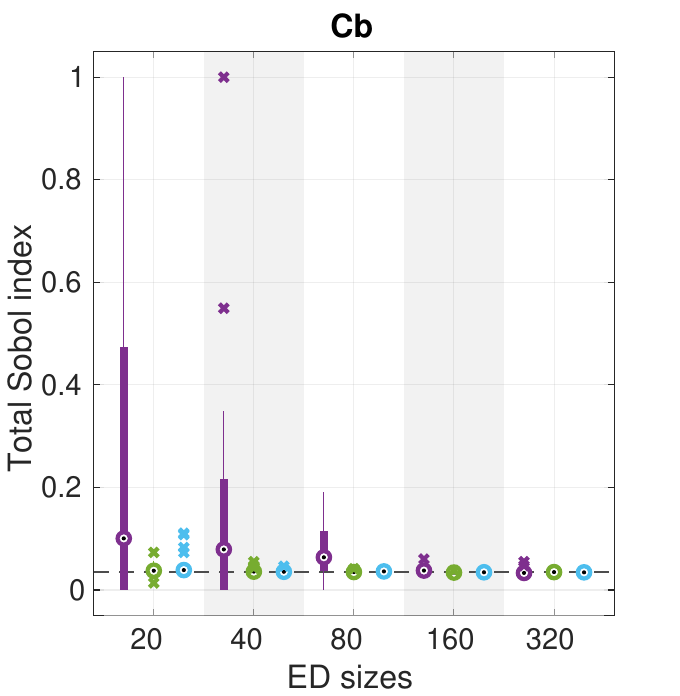}
  }
  \caption{Results for the flood cost model in the weighted case.}
  \label{fig:FloodModel:wPoinCE}
\end{figure}

\begin{figure}[h!]
    \centering  \includegraphics[width=0.7\linewidth]{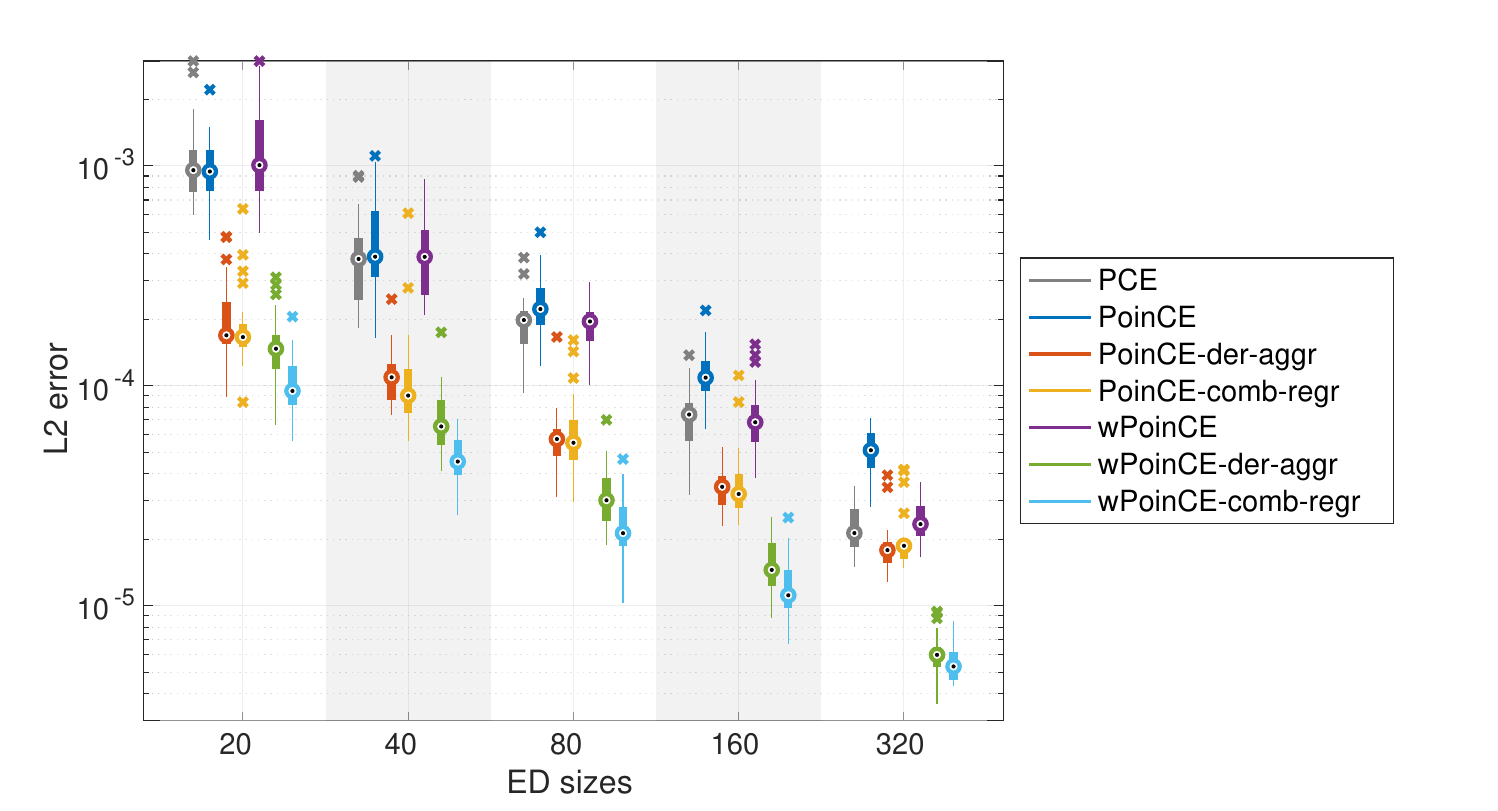}
    \caption{\rev{$\Ltwo$ error. Results for the flood model, both unweighted and weighted PoinCE}, \rev{as well as PCE as a comparison.}}
    \label{fig:FloodModel:PoinCE:L2}
\end{figure}

Finally, using weighted Poincar\'e offers a clear advantage in this model. Indeed, wPoinCE-comb-regr provides the most accurate total Sobol' indices estimations, as shown in Figures \ref{fig:FloodModel:PoinCE} and \ref{fig:FloodModel:wPoinCE}, and  yields the lowest $\Ltwo$ errors, displayed in Figure \ref{fig:FloodModel:PoinCE:L2}.

\rev{Figure \ref{fig:FloodModel:PoinCE:L2} also includes the $\Ltwo$ error for standard PCE without gradients. Comparing the convergence speed of PCE, (unweighted) PoinCE and wPoinCE, this plot illustrates the result from Section~\ref{sec:approx_properties}: weighted PoinCE can achieve the same convergence as polynomial approximations 
when the weight vanishes at the boundaries, which is fulfilled when using $w_\text{lin}$, but not for unweighted PoinCE ($w = 1$).
We do not show results for PCE with gradients, since in general the partial derivatives of a PCE basis are not orthogonal and therefore not suited to be directly included into a combined regression problem, as explained in Section~\ref{sec:sparse_regression}. Poincar\'e expansions uniquely combine fast convergence with the ability of including derivative information.
}

\section{Conclusion and perspectives}

We investigated the use of the Poincaré basis for gradient-enhanced sensitivity analysis with chaos expansions. We demonstrated that this orthonormal basis is stable under differentiation, meaning that its derivatives also form an orthogonal basis. This key property was used to obtain recovery guarantees in sparse regression problems \cite{Adcock_2019}, and has traditionally been leveraged in specific cases of gradient-enhanced polynomial chaos expansions.
\rev{We further show that under certain conditions, the Poincaré basis achieves the same convergence rate as the best polynomial approximation, for classes of smooth functions.}

We then introduced a general framework for gradient-enhanced global sensitivity analysis (GSA), combining efficient algorithms for sparse regression with construction of weights in Poincaré inequalities. Notably, the methodology does not require closed-form expressions for the Poincaré basis, making it applicable to a wide range of probability measures and significantly broadening the applicability of gradient-enhanced chaos expansions.

We assessed the performance of Poincaré chaos expansions on a toy model with interactions and on a challenging 8-dimensional flood model, where input variables follow various, mostly non-standard, distributions motivated by real-world considerations. Several conclusions emerged.
First, incorporating derivatives consistently improves both the construction of the expansion and the estimation of Sobol' indices.
Second, using the weight $w_\text{lin}$ in Poincaré inequalities is beneficial, as it leads to a basis that includes linear functions, which is particularly useful for capturing linear trends in models.
Third, when combining function and derivative evaluations to build the expansion, the aggregated estimator (from multiple single-output regressions) performs competitively with the estimator from multi-output regression, which is advantageous as the size of the experimental matrix increases.

Throughout the paper, we primarily focused on two weight choices for defining the Poincaré basis: the constant weight and the one corresponding to the existence of a linear basis function. However, more general weights can be considered and may further enhance GSA performance. For instance, \cite{HerediaWeightPoincare} proposes a weight derived from a monotonic approximation of the model’s main effect. Another promising direction is to extend the framework to dependent input variables. A first step could be to consider independent groups of variables, allowing dependencies within each group. For sufficiently small groups, numerical construction of the Poincaré basis remains feasible.

\rev{A further avenue of research would be to explore the conditions that ensure both the existence of the Poincaré basis and its desirable properties, such as derivative orthogonality and approximation. The current conditions are overly restrictive. For example, the sufficient condition that $1/p$ is integrable does not hold for Legendre polynomials, which are a specific instance of the Poincaré basis, 
yet these polynomials still satisfy the conclusions of Propositions~\ref{prop:PoincareBasis_and_Derivation} and \ref{prop:approximation_property}. This direction would require a more in-depth analysis of spectral theory.}


\section*{Acknowledgement and funding}
\rev{We are grateful to the associate editor and the two anonymous reviewers, whose feedback led to substantial improvements, particularly regarding the approximation properties.}
The authors declare that AI tools were used solely to refine the text in terms of spelling, grammar, and overall style.\\
This work has benefited from the AI Interdisciplinary Institute ANITI. ANITI is funded by the France 2030 program under the Grant agreement ANR-23-IACL-0002

\appendix

\section{Proofs} \label{sec:Appendix}
This appendix contains all the proofs that have been omitted in the main text. The proofs of Propositions \ref{prop:sufficient_conditions_for_existence} and \ref{prop:PoincareBasis_and_Derivation} rely on some regularity of functions in $\HoneW$, which is established in the following lemma.
\begin{lmm}
\label{lemma:absolutely_continuous_representative}
    Let $\mu \in \probClass, \weight \in \funClassWeight$ and assume that $1/(\weight\rho)\in \LoneLeb$. Then every function $f$ in $\HoneW$ admits a representative $\tilde{f}$ that is absolutely continuous, i.e. $f = \tilde{f}$ almost everywhere (with respect to the Lebesgue measure). For simplicity we still denote by $f$ (instead of $\tilde{f}$) this representative. Thus we have
    \begin{equation}
        \label{eq:FTC_absolutely_continuous} f(x)=f(y)+\int_a^b f'(z)\,dx\quad \mbox{for all }x,y\in (a,b).
    \end{equation}
\end{lmm}
\begin{proof}
 Let $f\in \HoneW$. Since $1/(w\rho)\in \LoneLeb$, we can apply the Cauchy-Schwarz inequality to show that $f'\in \LoneLeb$:
\[\int_a^b |f'|\leq \left(\int_a^b \weight \,(f')^2\rho \right)^{1/2}\left(\int_a^b \frac{1}{\weight\,\rho} \right)^{1/2}< +\infty.\]
Therefore for any $x_0\in (a,b)$ fixed, the function 
$x\mapsto F(x) := \int_{x_0}^x f'(z)\,dz$ 
is absolutely continuous and such that $F'=f'$ almost everywhere  \cite[Lemma 3.31]{absolutely_continuous_book}. Then there exists some constant $c\in \R$ such that $F+c$ is equal to $f$ almost everywhere. In particular \eqref{eq:FTC_absolutely_continuous} holds for this representative according to Theorem 3.30 in \cite{absolutely_continuous_book}.
\end{proof}
\begin{proof}[Proof of Proposition \ref{prop:sufficient_conditions_for_existence}] 
Consider the Sturm-Liouville eigenvalue problem
\begin{equation}
\label{eq:SturmLiouvilleProblemAnnexe}
    \left\{\begin{array}{ll}
    \displaystyle L_w(f)=\frac{1}{\rho}(p f')'=-\eigval f\rho,\quad \mbox{on }(a,b),\\
     (p f')(a)=(p f')(b)=0.
\end{array}
\right.
\end{equation}
with $\p = \weight \, \rho$. 
Assume Condition (i): $1/p \in \LoneLeb$. In this case, Problem \eqref{eq:SturmLiouvilleProblemAnnexe} is said to be regular and there exists a countable family of eigenfunctions $(\psi_\indBasis)_{\indBasis \in \N}$, with respective eigenvalues $0 = \lambda_0 < \lambda_1 < \dots < \lambda_\indBasis < \dots $,  \cite[Theorem 4.3.1, (6)]{Zettl}. The functions $\psi_\indBasis$ and $\p\,\psi_\indBasis'$ are absolutely continuous in $(a,b)$ (this is included in the definition of an eigenfunction in \cite[Definition 2.2.1]{Zettl}) and the family $(\psi_\indBasis)_{\indBasis \in \N}$ forms a basis in $\Ltwo$, according to Theorem 4.11.1 in the same reference. Notice that, since the eigenvalues are simple, the eigenfunctions are uniquely defined, up to normalization and sign change. 

Since integration by parts applies to absolutely continuous functions \cite[Corollary 3.37]{absolutely_continuous_book}, we can deduce that each $\psi_\indBasis'$ belongs to $\LtwoW$. Indeed, since each $\psi_\indBasis$ satisfies the boundary conditions in \eqref{eq:SturmLiouvilleProblemAnnexe} we have
\[\int_a^b w(\psi_\indBasis')^2\rho=-\int_a^b (\p\,\psi_\indBasis')'\, \psi_\indBasis=\int_a^b (-L_w(\psi_\indBasis))\,\psi_\indBasis\,\rho=\lambda_\indBasis\int_a^b (\psi_\indBasis)^2\rho <+\infty.\]
Finally, let $g \in \HoneW$. Due to Lemma \ref{lemma:absolutely_continuous_representative} we can assume that $g$ is absolutely continuous. Thus, we can use a similar  integration by parts to show that
\[\lambda_\indBasis\langle \psi_\indBasis,g\rangle=\langle \psi_\indBasis',g'\rangle_w,
\]
which proves the existence of the Poincar\'e basis.\\

Now, suppose instead that condition $(ii)$ holds: that is, that the anti-derivatives of $1/\p$ belong to $\Ltwo$. We fix one of them and denote it by $R$. Observe that when we fix $\eigval=0$ and ignore the boundary conditions in \eqref{eq:SturmLiouvilleProblemAnnexe}, all the solutions consist of linear combinations of the constant function 1 and $R$. \rev{Then, each solution belongs to $\Ltwo$ and has at most one zero in $(a,b)$ since $R$ is strictly monotonic ($R'=1/p>0$ in $(a,b)$).
This implies in particular that} 
the endpoints $a,b$ are limit-circle non-oscillatory, according to Sturm-Liouville formalism 
\cite[Definition 7.3.1]{Zettl}, and Theorem 10.12.1(4) in
this reference guarantees the
existence of $(\psi_\indBasis)_\indBasis$, the countable family of eigenfunctions of \ref{eq:SturmLiouvilleProblemAnnexe}, such that $\psi_\indBasis$ and $p\,\psi'_\indBasis$ are absolutely continuous. 

However this result does not mention the fact that $(\psi_\indBasis)_{\indBasis\in \N}$ forms a basis in $\Ltwo$. This can be deduced from its proof in \cite{Zettl_singular}. Indeed, reading Section 5 in this reference, there exists some positive function $v$ in $(a,b)$ such that the functions $\phi_{\indBasis}=\psi_\indBasis/v$ are solutions of a regular Sturm-Liouville problem associated to the operator 
${L_{w,v} f=\frac{1}{v^2\rho}(v^2\p f')'+\frac{1}{v\,\rho}\left(p \,v'\right)f}$. In particular $(\phi_{\indBasis})_{j\in \N}$ is a basis in $\xLtwo(\mu,v^2)$, which is equivalent to say that $(\psi_{\indBasis})_{j\in \N}$ is a basis in $\Ltwo$.  The remainder of the proof proceeds exactly as in the regular case.
\end{proof}

\begin{proof}[Proof of Proposition \ref{prop:PoincareBasis_and_Derivation}]
First notice that
the assumptions on $\mu, \weight$
ensure the existence of the Poincar\'e basis, by Proposition~\ref{prop:sufficient_conditions_for_existence}, ($i$).

Let us start by the direct sense (1).
The fact that $(\eigfun'_\indBasis)$ is an orthogonal system of $\LtwoW$ directly comes from \eqref{eq:PoincareWeakProblem}. Indeed, by choosing $f =\eigfun_\indBasis$, $g = \eigfun_m$ and $\eigval = \eigval_\indBasis$, we obtain 
$$ \LtwoDotW{\eigfun_\indBasis}{\eigfun_m} = \delta_{\indBasis, m} \eigval_\indBasis.$$
Consequently, we obtain 
$$ \LtwoDotH{\eigfun_\indBasis}{\eigfun_m} = \delta_{\indBasis, m} (1+\eigval_\indBasis)$$
which also proves that $(\eigfun_\indBasis)$ is an orthogonal system of $\HoneW$.\\
It remains to show that these two systems are complete.\\
Let us prove first that the space spanned by $(\eigfun_\indBasis)$ is dense in $\HoneW$, by proving that its orthogonal is the null space. Let $f \in \HoneW$ such that $\LtwoDotH{f}{\eigfun_\indBasis}=0$ for all $\indBasis \in \N$. Then, by \eqref{eq:PoincareWeakProblem}, this implies $(1+\eigval_\indBasis) \LtwoDot{f}{\eigfun_\indBasis}=0$. Thus $f=0$ as $1+\eigval_\indBasis >0$ and $(\eigfun_\indBasis)$ is a basis of $\Ltwo$.

Now, to show that the span of $(\varphi'_\indBasis)$ is dense in $\LtwoW$ let us further assume that $1/\p \in \LoneLeb$.
 Let $f \in \LtwoW$ be a function such that $\LtwoDotW{f}{\eigfun'_\indBasis} = 0$ for all $\indBasis \geq 1$. 
Due to our assumptions, the anti-\rev{derivative} $x\mapsto g(x)= \int_a^x f(y) \,dy$ belongs to $\HoneW$. Indeed, first we have $g'=f\in \LtwoW$. Second, we can prove that $g$ is bounded and thus belongs to $\Ltwo$ by using the following Cauchy-Schwarz inequality:
\[ \vert g \vert \leq \int_a^b | f | \ \leq \left(\int_a^b \weight\, f^2 \rho\ \right)^{1/2}\left( \int_a^b \frac{1}{\weight \, \rho} \right)^{1/2}<\infty.\]

Hence the condition $\LtwoDotW{f}{\eigfun'_\indBasis} = 0$ is rewritten $\LtwoDotW{g'}{\eigfun'_\indBasis} = 0$ and by \eqref{eq:PoincareWeakProblem}, $\LtwoDot{g}{\eigfun_\indBasis} = 0$ for $\indBasis \geq 1$. We conclude that $g$ is proportional to $\eigfun_0$, which is a constant function. Thus $f=g'=0$.\\

Finally, still assuming that $ 1/\p\ \in \LoneLeb$, 
let us consider an orthonormal basis $(\psi_\indBasis)$ of $\Ltwo$ contained in $\HoneW$, with $\psi_0\equiv 1$ and such that $(\psi'_\indBasis)$ is an orthogonal basis in $\LtwoW$. Let us prove that it coincides with the Poincaré basis.

Recall that we can replace any function of $\HoneW$ by its continuous representative (see Lemma \ref{lemma:absolutely_continuous_representative}).
Now, for every $f\in \HoneW$ we have the expansions
\begin{equation}
\label{eq:bases_expansions_in_proposition_poincare_derivative_basis}
    f=\sum_{\indBasis\geq 0}\langle f,\psi_\indBasis\rangle\psi_\indBasis\quad \mbox{and}\quad f'=\sum_{\indBasis\geq 1}\frac{\langle f',\psi'_\indBasis\rangle_w}{\norme{\psi_\indBasis'}{\weight}^2} \,\psi_\indBasis'.
\end{equation}
Moreover, given $x_0\in (a,b)$ we have the following alternative expression for $f$:
\begin{equation}
\label{eq:identity_for_f_in_proposition_derivative_basis}
    f(x)=f(x_0)+\sum_{\indBasis\geq 1}\frac{\langle f',\psi'_\indBasis\rangle_w}{\norme{\psi_\indBasis'}{\weight}^2} \left(\psi_\indBasis(x)-\psi_\indBasis(x_0)\right),\quad \mbox{for all } x\in (a,b).
\end{equation}
Indeed, given $N\geq 1$, \eqref{eq:FTC_absolutely_continuous} gives
\begin{align*}
    &\left|f(x)-f(x_0)-\sum_{\indBasis=1}^N \frac{\langle f',\psi'_\indBasis\rangle_w}{\norme{\psi_\indBasis'}{\weight}^2}(\psi_\indBasis(x)-\psi_\indBasis(x_0))\right| \\
    &=\left|\int_{x_0}^x \left(f'(z)-\sum_{\indBasis=1}^N\frac{\langle f',\psi'_\indBasis\rangle_w}{\norme{\psi_\indBasis'}{\weight}^2} \psi_\indBasis'(z)\right)\,dz\right|
    \leq \norme{f'-\sum_{\indBasis=1}^N\frac{\langle f',\psi'_\indBasis\rangle_w}{\norme{\psi_\indBasis'}{\weight}^2} \psi_\indBasis'}{\weight}\left(\int_a^b \frac{1}{\weight\rho}\right)^{1/2}.
\end{align*}
The right-hand side converges to 0 as $N\rightarrow \infty$ and thus we obtain \eqref{eq:identity_for_f_in_proposition_derivative_basis}.
Now, since $\psi_0\equiv 1$ and $(\psi_\indBasis)$ is an orthonormal basis in $\Ltwo$, after multiplying each side by $\psi_m$ ($m\geq 1$) and integrating with respect to $\mu$ we obtain
\[\lambda_m \langle f,\psi_m\rangle=\langle f', \psi'_m\rangle_w,\quad\mbox{for all }f\in \HoneW,\]
with $\lambda_m=\norme{\psi_m'}{\weight}^2>0$. This shows that $(\psi_m)_{m \in \N}$ is the Poincar\'e basis.
\end{proof}

\begin{proof}[Proof of Proposition~\ref{prop:DGSMfromPoinCE}]
    From the proof of Proposition~\ref{prop:PoincareBasis_and_Derivation}, we can write
    $$ \partderiv[\cm]{\indVar}(x) = \sum_{\alp \in \N^d, \alpha_\indVar \geq 1} c_\alp 
    \partderiv[\psi_\alp]{\indVar} $$
    Now, $\partderiv[\psi_\alp]{\indVar}(x) = \psi_{\indVar, \alpha_\indVar}'(x_\indVar) \prod_{\ell \neq \indVar} \psi_{\ell, \alpha_\ell}(x_\ell) $. As the derivatives of the Poincaré basis functions form an orthogonal basis, we can use a Parseval identity, given $X_\ell$ with $\ell \neq \indVar$, to obtain
    \begin{equation*} \label{eq:DGSMaux}
\Esp_{X_\indVar} \left[ w_\indVar(X_\indVar) \left(  \partderiv[\cm]{\indVar}(X) \right)^2 \right] = \sum_{\alp \in \N^d, \alpha_\indVar \geq 1} c_\alp^2 \Esp\left[w_\indVar(X_\indVar) \psi_{\indVar, \alpha_\indVar}'(X_\indVar)^2\right] \prod_{\ell \neq \indVar} \psi_{\ell, \alpha_\ell}(X_\ell)^2        
    \end{equation*}
Notice that $\Esp\left[w_\indVar(X_\indVar)\psi_{\indVar, \alpha_\indVar}'(X_\indVar)^2\right] = \LtwoNormW{\psi_{\indVar, \alpha_\indVar}}^2= \eigval_{\indVar, \alpha_\indVar}$, and for all $\ell = 1, \dots, d$,
$\Esp\left[\psi_{\ell, \alpha_\ell}(X_\ell)^2\right] = \LtwoNorm{\psi_{\ell, \alpha_\ell}}^2= 1$. Then, by integrating the expression above over $X_\ell$ with $\ell \neq \indVar$, we get
$$ \nu_\indVar = \sum_{\alp \in \N^d, \alpha_\indVar \geq 1} \eigval_{\indVar, \alpha_\indVar} c_\alp^2$$
To obtain the inequality between $\nu_\indVar$ and $S_\indVar^\tot$, it is sufficient to use the inequality $\eigval_{\indVar, \alpha_\indVar} \geq \eigval_{\indVar, 1}$, valid for all $\alpha_\indVar \geq 1$:
$$ \nu_\indVar \geq \rev{\eigval_{\indVar, 1}} \times \sum_{\alp \in \N^d, \alpha_\indVar \geq 1} c_\alp^2 = \frac{1}{C_P(\mu_\indVar, \weight_\indVar)} S_\indVar^\tot \, \Var(\cm(X)).$$
\end{proof}

\rev{
\begin{proof}[Proof of Proposition 
\ref{prop:approximation_property}]
The scheme of the proof is inspired by  \cite[\S 5.2.2]{canuto2006spectral}. We adapt it to our assumptions.
\begin{enumerate}
    \item [(i)] First, note that due to the hypotheses on $w$, $\rho$, and $f$, the function $(-L_w)^k f$ belongs to $\Ltwo$ for each $k$ such that $1 \leq k \leq 2m$. Indeed, the operator $L_w$ can be written as
    $$ L_wf = \frac{1}{\rho}(w\,\rho\,f')' = wf'' + \tau f', $$
    where we recall that the functions $w$ and $\tau = (w\,\rho)'/\rho$ are of class $\Cm{2k-2}$. Thus starting from $f \in \Hm{2k}$ we get $L_wf \in \Hm{2k-2}$. Iterating $k$ times, we obtain that $(L_w)^k f \in \Ltwo$.\\
    Now, by derivating $k$ times, we obtain the expansion
    \[ (-L_w)^k f = \sum_{i=1}^{2k} \Phi_{k,i} f^{(i)}, \]
    where each $\Phi_{k,i}$
    is a linear function of $w, \tau$ and their
    derivatives up to order $2k-2$.
    With our assumptions, each $\Phi_{k,i}$ is thus bounded on $[a,b]$.
    To bound the $\Ltwo$ norm, we first apply the Cauchy-Schwarz inequality for sums:
    \[ \left((-L_w)^k f\right)^2 \leq \left( \sum_{i=1}^{2k} |\Phi_{k,i}\, f^{(i)}| \right)^2 \leq \left( \sum_{i=1}^{2k} \Phi_{k,i}^2 \right) \left( \sum_{i=1}^{2k} \left(f^{(i)}\right)^2 \right). \]
    Integrating with respect to $\mu$ gives 
    \begin{equation} \label{eq:LiterateNormBound}
    \norme{(-L_w)^k f}{}\leq C_k \norme{f}{\Hm{2k}},
    \end{equation} 
    with 
    $C_k = \left\| \sqrt{\sum_{i=1}^{2k} \Phi_{k,i}^2} \right\|_\infty$, 
    which shows the continuity of $(L_w)^k: \Hm{2k} \to \Ltwo$.\\
    Let us now derive \eqref{eq:coefFormulaWithLiterate}.
   By \eqref{eq:PoincareWeakProblem}, we have
    $$ \hat{f}_\indBasis = \langle f, \eigfun_\indBasis \rangle 
    = \frac{1}{\lambda_\indBasis} \langle f', \eigfun'_\indBasis \rangle_w 
    = \frac{1}{\lambda_\indBasis} \int_a^b   (w r f') \eigfun'_\indBasis  dx. $$
    Since $f \in \Hm{2k}$ and the operator coefficients are smooth, we can integrate by part. This gives
    $$ \hat{f}_\indBasis = -\frac{1}{\lambda_\indBasis} \int_a^b (w \rho f')' \eigfun_\indBasis \ dx 
    =  \frac{1}{\lambda_\indBasis} \langle -L_w f, \eigfun_\indBasis \rangle $$
    where the boundary term again vanishes due to the condition $p(a)=p(b)=0$ (with $p=wr$). 
    By (i), $L_w f \in \Hm{2k-2}$. 
    Thus we can iterate this process $k$ times, leading to \eqref{eq:coefFormulaWithLiterate}. 
    \item[(ii)]  
    This claim follows from (i) and known convergence properties of the eigenvalues $(\lambda_\indBasis)_{\indBasis\geq 0}$. 
    Indeed, \cite[Theorem 4.3.1, (6)]{Zettl} tells that
    \begin{equation}
    \label{eq:eigenvalue_convergence}
        \lim_{\indBasis\rightarrow \infty}\frac{\lambda_\indBasis}{\indBasis^2}=\pi^2 \left(\int_a^b 
    \frac{1}{\sqrt{w}}
    \right)^{-2}
    \end{equation}
    Notice that $\displaystyle \frac{1}{\sqrt{w}}= \frac{1}{\sqrt{p}} \sqrt{r}$ belongs to $\xLone(a,b)$ as a product of two elements of $\xLtwo(a,b)$.\\
    Therefore, combining with \eqref{eq:coefFormulaWithLiterate} and using Cauchy-Schwarz inequality, one can find $C>0$ such that for all $k$ such that $2k\leq m$ and all $j \geq 1$,
    \[|\hat{f}_\indBasis| \leq 
    \frac{C}{\indBasis^{2k}}\norme{(-L_w)^k f}{}\]
\clearpage

    \item[(iii)] From Parseval inequality, \rev{the fact that the sequence $(\eigval_{j})_{j \geq 1}$ is increasing}, and \eqref{eq:coefFormulaWithLiterate} we obtain
    \begin{eqnarray*}
    \norme{f-P_Nf}{}^2
    =\sum_{\indBasis=N+1}^\infty | \hat{f}_\indBasis |^2
    &=&  \sum_{\indBasis=N+1}^\infty \frac{1}{\eigval_\indBasis^{2k}} |\langle (-L_w)^k f,\eigfun_\indBasis  \rangle|^2\\
    &\leq& \frac{1}{\eigval_{N+1}^{2k}} \sum_{\indBasis=N+1}^\infty  |\langle (-L_w)^k f,\eigfun_\indBasis  \rangle|^2
    \leq \frac{1}{\eigval_{N+1}^{2k}}\norme{(-L_w)^kf}{}^2.
    \end{eqnarray*}
    The result follows by
    combining with \eqref{eq:eigenvalue_convergence} and the estimate \eqref{eq:LiterateNormBound}. 
\end{enumerate}
\end{proof}
}

\rev{
\section{Gradient-informed computation of the expansion coefficients}
\label{sec:app:regression}
}

In practice, truncated chaos expansions are computed from a finite set of model evaluations at specified points from the input domain called the \emph{experimental design}.
In this paper, we always sample the experimental design from the joint pdf of the input random variables. Other choices are possible, for example the so-called coherence-optimal sampling \cite{Hampton2015,LuethenSIAMJUQ2020}, also analyzed by \cite{Adcock_2019}.

Let $N$ denote the number of model evaluations and $P$ the number of coefficients in the truncated expansion.
From the various available methodologies, including collocation and (sparse) quadrature, we choose regression-based methods which are sample-efficient and stable.

We use \emph{sparse regression} to compute the chaos coefficients, i.e., ordinary least squares regression with a regularization term which enforces sparsity in the chaos coefficients (usually $\ell^1$ minimization). 
To this aim, we assemble the \emph{regression matrix} $\ve \Psi \in \R^{N \times P}$ consisting of evaluations of the basis functions at the experimental design points:
\[
    \Psi_{\indDoE, \indBasis} = \psi_{\indBasis}(\ve x^{(\indDoE)}), \indDoE = 1 \enum N, \indBasis = 1 \enum P
\]
and the vector of model evaluations
\[
    \ve y = (\cm(\ve x^{(1)}) \enum \cm(\ve x^{(N)}))^T.
\]
Denoting by $\ve c = (c_{\ve\alpha_1} \enum c_{\ve\alpha_P})^T$ the vector of chaos coefficients, we are looking for a $\ve c$ that is sparse while fulfilling $\ve \Psi \ve c \approx \ve y$.

There is a wide variety of sparse regression methods, ranging from greedy stepwise algorithms to Bayesian techniques as reviewed in \cite{LuethenSIAMJUQ2020, LuethenIJUQ2022}. We apply least-angle regression (LARS) model selection achieved by using the leave-one-out error \cite{UQdoc_21_104}. 
Note that Adcock and Sui \cite{Adcock_2019} use the SPGL1 solver, which in our benchmark \cite{LuethenSIAMJUQ2020} however did not perform as well as most other tested sparse regression solvers.

\rev{
In the following sections we provide more details on two different ways of incorporating derivative information into the computation process.
}


\rev{
\subsection{Aggregating derivative expansions}
\label{app:averaged_deriv_expansions}
}

In practice, derivative evaluations can be utilized in different ways to compute the chaos coefficients. In \cite{Poincare_Chaos_Sparse_2023} we presented one possible approach, which we recall here. Let 	
\begin{equation}
	 \cm(\ve x) 
	\ \approx \ 
    \tilde\cm(\ve x) 
	= \sum_{\alp\in \ca} c_\alp \psi_\alp(\ve x)
\end{equation}
denote the finite chaos expansion. Then by partial differentiation we get
\begin{equation}
\label{eq:regression_derivative_evaluations}
    \partderiv[\cm]{\indVar} (\ve x)
	\ \approx \ 
    \partderiv[\tilde\cm]{\indVar} (\ve x) 
	= \sum_{\alp\in \ca: \alpha_\indVar > 0} c_\alp \partderiv[\psi_\alp]{\indVar}(\ve x).
\end{equation}
Since basis functions which are constant in $x_\indVar$ vanish by differentiation with respect to $x_\indVar$, the $\indVar$-th derivative expansion only involves the subset of terms $\{\alp \in \ca: \alpha_\indVar > 0\}$, and the remaining coefficients cannot be estimated from the derivative data.

Including the weight function into \eqref{eq:regression_derivative_evaluations}, we can write 
\begin{equation}
\label{eq:regression_derivative_evaluations_weighted}
    \sqrt{w_\indVar(\ve x)}\partderiv[\cm]{\indVar} (\ve x)
	\ \approx \ 
    \sqrt{w_\indVar(\ve x)}\partderiv[\tilde\cm]{\indVar} (\ve x) 
	= \sum_{\alp\in \ca: \alpha_\indVar > 0} c_\alp \sqrt{w_\indVar(\ve x)} \partderiv[\psi_\alp]{\indVar} (\ve x).
\end{equation}
As explained in Section~\ref{sec:poincare_basis}, the derivatives form again an orthogonal basis in $\LtwoW$; alternatively, we can see that $\{\sqrt{w_\indVar(\ve x)} \partderiv[\psi_\alp]{\indVar} (\ve x)\}$ forms an orthogonal basis in $\Ltwo$. 
In our aggregation approach, sparse regression (Section~\ref{sec:sparse_regression}) is applied as usual to compute the coefficients for each of the $d$ derivative expansions.
Note that this is a unique advantage of using Poincar\'e basis functions: we cannot do the same with orthogonal polynomials, as their derivatives do in general not form an orthogonal basis in any $\xLtwo$ space (except for Hermite, Laguerre and Jacobi polynomials).

When applying sparse regression to the $d$ derivative expansions in Eq.~\eqref{eq:regression_derivative_evaluations}, the coefficient estimates computed from different expansions do in general not coincide due to the finite size of the data set. 
Let $\hat{\ve c}^{\delta, k}$ denote the vector of coefficient estimates of the $k$-th derivative expansion.
Our aggregation approach determines new coefficient estimates $\hat{c}^{\text{aggr}}_{\alp_1} \enum \hat{c}^{\text{aggr}}_{\alp_P}$ as follows:
\begin{itemize}
    \item Coefficients of univariate terms are estimated by one single derivative expansion and this value is used.
    \item Coefficients corresponding to bivariate terms are estimated by two derivative expansions. They are averaged to yield the final estimate.
    \item More generally, for $m = 1 \enum d$, coefficients corresponding to rank-$m$ terms are estimated by $m$ derivative expansions. They are averaged to yield the final estimate.
    \rev{Here, the term ``rank'' denotes the number of active input variables, i.e., the number of variables in which this polynomial is not constant.}
    \item The constant coefficient cannot be estimated by any of the derivative expansions. Its computation is explained in Eqs.~\eqref{eq:aggregated_constant_coeff1}-\eqref{eq:aggregated_constant_coeff2}.
\end{itemize}
To put this idea in mathematical terms, let $\hat{\ve c}^{\partial, \indVar}$ denote the vector of chaos coefficients computed from the $\indVar$-th derivative expansion.
The final coefficient estimate, averaged over all contributing derivative expansion results, is computed as follows:
\begin{equation}
\hat{c}_\alp^{\text{aggr}} = \frac{1}{|\{\indVar \in \{1 \enum d\}:\alpha_\indVar \geq 1\}|}\sum_{\indVar: \alpha_\indVar \geq 1} \hat{c}_\alp^{\, \partial, \indVar} \quad \text{ for each } \alp \in \ca \setminus \ve{0}.
\end{equation}
The set of coefficients $\ca \setminus \ve{0}$ is sufficient to compute Sobol' and DGSM indices, as visible from Eqs.\eqref{eq:totalSobolFromCE} and \eqref{eq:DGSMfromPoinCE}.
If the constant coefficient is needed, for example when the chaos expansion should be used as a surrogate model, it needs to be computed from the model evaluations only, for example, using ordinary least-squares on the residual $\ve{y}_\text{res}$:
\begin{align}
\label{eq:aggregated_constant_coeff1}
\ve{y}_\text{res} &= \ve{y} - \ve\Psi \begin{pmatrix} 0 \\ \hat{c}_{\alp_1}^{\, \partial, \text{aggr}} \\ \vdots \\ \hat{c}_{\alp_{P-1}}^{\, \partial, \text{aggr}} \end{pmatrix}, \\
\label{eq:aggregated_constant_coeff2}
\hat{c}_{\ve{0}}^{\, \partial, \text{aggr}} &= \frac{1}{N} \sum_{\indDoE=1}^{N} \ve{y}_\text{res}{(\indDoE)}.
\end{align}

\rev{
The main advantage of this approach is that 
it requires the solution of several small regression problems instead of assembling, weighting and normalizing one large problem as described in Section~\ref{app:gradient_enhanced_regression}.}


\rev{
\subsection{Gradient-enhanced regression (combined regression)}
\label{app:gradient_enhanced_regression}
}
Another way to utilize derivative evaluations is to set up a regression problem that combines all available data: model evaluations as well as partial derivative evaluations. 
This idea was already mentioned in \cite{Jakeman2015} and developed in generality by \cite{Adcock_2019}.

Recall that the original regression problem (based on model evaluations only) is given by
\begin{equation}
	 \cm(\ve x) 
	\rev{~\approx~} \sum_{\alp\in \ca} c_\alp \psi_\alp(\ve x) \qquad \Rightarrow \qquad \ve \Psi \ve c \rev{~\approx~} \ve y,
\end{equation}
with the regression matrix $\{\Psi_{\indDoE, \indBasis} = \psi_\indBasis(\ve x^{(\indDoE)}), \, \indDoE = 1, \cdots,N; \; \indBasis=1, \cdots, P \equiv \textrm{card} \mathcal{A}\}$ and model responses $y_i = \cm(\ve x^{(\indDoE)})$.

Consider first the case of Poincar\'e expansions without weights ($w_\indVar \equiv 1$ for all $\indVar\in\{1,\dots,d\}$).
The gradient-enhanced regression problem, which we call \emph{combined regression} problem, is given by
\begin{equation}
    \begin{bmatrix}
        \ve \Psi \\
        \ve \Psi_{\partial, 1} \\
        \vdots \\
        \ve \Psi_{\partial, d}
    \end{bmatrix}
    \ve c 
     \rev{~\approx~} 
     \begin{bmatrix}
         \ve y\\
         \ve y_{\partial, 1}\\
         \vdots \\
         \ve y_{\partial, d}
     \end{bmatrix}
\end{equation}
with regression matrices $(\ve \Psi_{\partial, \indVar})_{\indDoE, \indBasis} = \partderiv[\psi_\indBasis]{\indVar} (\ve x^{(\indDoE)})$ and model derivatives $(\ve y_{\partial, \indVar})_i = \partderiv[\cm]{\indVar} (\ve x^{(\indDoE)})$.

In case of weighted Poincar\'e expansions, $w_\indVar \neq 1$, the regression problem needs to be modified by pre-multiplication with diagonal scaling matrices $\ve T_\indVar$:
\begin{equation}
    \begin{bmatrix}
        \ve \Psi \\
        \ve T_1 \ve \Psi_{\partial, 1} \\
        \vdots \\
        \ve T_d \ve \Psi_{\partial, d}
    \end{bmatrix}
    \ve c 
     \rev{~\approx~}
     \begin{bmatrix}
         \ve y\\
         \ve T_1 \ve y_{\partial, 1}\\
         \vdots \\
         \ve T_d \ve y_{\partial, d}
     \end{bmatrix}
\end{equation}
with entries $(T_\indVar)_{\indDoE,\indDoE} = \sqrt{w_\indVar(x_\indVar^{(\indDoE)})}$,
where $x_\indVar^{(\indDoE)}$ is the $\indVar$-th entry of the $\indDoE$-th sample point.
This procedure is sometimes referred to as \emph{preconditioning} \cite{Peng2016,Guo2018}.
Note that $\ve T_\indVar$ becomes slightly more involved if the sampling density does not coincide with the joint density $\rho$ of the input random variables \cite{Adcock_2019}.

Note that this preconditioning is particularly simple: it only involves diagonal matrices. This is specific to the Poincar\'e basis, because its derivatives form an orthogonal basis, ensuring that all the regression matrices $\Psi_{\partial, \indVar}$ are diagonal in expectation (when the sampling density is equal to $\rho$).

Finally, in both cases each column of the regression matrix is divided by its norm to ensure that the columns of the resulting regression matrix are orthonormal in expectation (in $\HoneW$). With Proposition~\ref{prop:PoincareBasis_and_Derivation}, the norm of the column corresponding to multi-index $\alp$ is given by $\sqrt{1 + \sum_{\indVar=1}^d \eigval_{\indVar, \alpha_\indVar}}$.
This preconditioned regression problem comes with recovery guarantees when the coefficient vector $\ve c$ is sparse, see \cite[Section 4]{Adcock_2019}.

\begin{rmrk}
Notice that the framework in \cite{Adcock_2019} is slightly more restrictive than ours, since they consider the case of eigenfunctions of Sturm-Liouville problems \ref{eq:SturmLiouvilleProblem} with more restrictive boundary conditions $\chi(a) = \chi(b) = 0$. However, inspecting their proofs reveals that these boundary conditions are used to deduce that the eigenfunctions derivatives form an orthogonal basis, as in our Proposition~\ref{prop:PoincareBasis_and_Derivation}, on which the other demonstrations rely.
Therefore, the recovery guarantees are also valid for all Poincaré bases satisfying the assumptions of Proposition~\ref{prop:PoincareBasis_and_Derivation}.
\end{rmrk}

\bibliographystyle{abbrv} 
	\bibliography{bibli}

\end{document}